\documentclass[oneside,english]{acmart}
\usepackage[T1]{fontenc}
\usepackage[latin9]{inputenc}
\usepackage{array}
\usepackage{float}
\usepackage{units}
\usepackage{url}
\usepackage{multirow}
\usepackage{amsmath}
\usepackage{graphicx}

\makeatletter

\providecommand{\tabularnewline}{\\}
\floatstyle{ruled}
\newfloat{algorithm}{tbp}{loa}
\providecommand{\algorithmname}{Algorithm}
\floatname{algorithm}{\protect\algorithmname}

\theoremstyle{definition}
\newtheorem{defn}{\protect\definitionname}
\theoremstyle{plain}
\newtheorem{prop}{\protect\propositionname}
\theoremstyle{plain}
\newtheorem{thm}{\protect\theoremname}
\theoremstyle{definition}
 \newtheorem{example}{\protect\examplename}
\theoremstyle{remark}
\newtheorem{rem}{\protect\remarkname}

\AtBeginDocument{%
  \providecommand\BibTeX{{%
    \normalfont B\kern-0.5em{\scshape i\kern-0.25em b}\kern-0.8em\TeX}}}

\setcopyright{acmcopyright}

\acmJournal{TOMS}

\author{Qiao Chen}
\email{qiao.chen@stonybrook.edu}
\author{Xiangmin Jiao}
\email{xiangmin.jiao@stonybrook.edu}
\orcid{0000-0002-7111-9813}
\affiliation{%
  \institution{Stony Brook University}
  \department[0]{Department of Applied Mathematics \& Statistics and Institute for Advanced Computational Science}
  \city{Stony Brook}
  \state{New York}
  \country{USA}
  \postcode{11794}
}


\usepackage{algorithmic}
\usepackage{listings}

\def\vec#1{\boldsymbol{#1}}

\makeatother

\usepackage{babel}
\usepackage{listings}
\providecommand{\definitionname}{Definition}
\providecommand{\examplename}{Example}
\providecommand{\propositionname}{Proposition}
\providecommand{\remarkname}{Remark}
\providecommand{\theoremname}{Theorem}

\begin{document}
\title{HIFIR: Hybrid Incomplete Factorization with Iterative Refinement for
Preconditioning Ill-conditioned and Singular Systems}
\begin{abstract}
We introduce a software package called \emph{HIFIR} for preconditioning
sparse, unsymmetric, ill-conditioned, and potentially singular systems.
HIFIR computes a \emph{hybrid incomplete factorization}, which combines
multilevel incomplete LU factorization with a truncated, rank-revealing
QR factorization on the final Schur complement. This novel hybridization
is based on the new theory of \emph{$\epsilon$-accurate approximate
generalized inverse}. It enables near-optimal preconditioners for
consistent systems and enables flexible GMRES to solve inconsistent
systems when coupled with iterative refinement. In this paper, we
focus on some practical algorithmic and software issues of HIFIR.
In particular, we introduce a new inverse-based rook pivoting into
ILU, which improves the robustness and the overall efficiency for
some ill-conditioned systems by significantly reducing the size of
the final Schur complement for some systems. We also describe the
software design of HIFIR in terms of its efficient data structures
for supporting rook pivoting in a multilevel setting, its template-based
generic programming interfaces for mixed-precision real and complex
values in C++, and its user-friendly high-level interfaces in MATLAB
and Python. We demonstrate the effectiveness of HIFIR for ill-conditioned
or singular systems arising from several applications, including the
Helmholtz equation, linear elasticity, stationary incompressible Navier--Stokes
equations, and time-dependent advection-diffusion equation.
\end{abstract}

\maketitle

\renewcommand{\shorttitle}{HIFIR: Hybrid Incomplete Factorizations with Iterative Refinement}

\section{Introduction}

We consider the problem of preconditioning an iterative solver for
a linear system,
\begin{equation}
\boldsymbol{A}\boldsymbol{x}=\boldsymbol{b},\label{eq:linear-sys}
\end{equation}
where $\boldsymbol{A}\in\mathbb{C}^{n\times n}$ is sparse and potentially
singular, $\boldsymbol{x}\in\mathbb{C}^{n}$, and $\boldsymbol{b}\in\mathbb{C}^{n}$,
or as in many applications, $\boldsymbol{A}\in\mathbb{R}^{n\times n}$,
$\boldsymbol{x}\in\mathbb{R}^{n}$, and $\boldsymbol{b}\in\mathbb{R}^{n}$.
For generality, we will assume complex-valued systems in our discussions.
In general, (\ref{eq:linear-sys}) is \emph{inconsistent} in that
$\left\Vert \boldsymbol{b}-\boldsymbol{A}\boldsymbol{A}^{+}\boldsymbol{b}\right\Vert \gg\epsilon_{\text{mach}}\left\Vert \boldsymbol{b}\right\Vert $,
where $\boldsymbol{A}^{+}$ denotes the Moore--Penrose pseudoinverse
of $\boldsymbol{A}$ and $\epsilon_{\text{mach}}$ denotes the machine
epsilon of a given floating-point number system. In this case, we
must seek a \emph{least-squares solution} of (\ref{eq:linear-sys}),
i.e.,
\begin{equation}
\boldsymbol{x}_{\text{LS}}=\arg\min_{\boldsymbol{x}}\left\Vert \boldsymbol{b}-\boldsymbol{A}\boldsymbol{x}\right\Vert _{2},\label{eq:ls-sol}
\end{equation}
or often preferably the \emph{pseudoinverse} \emph{solution} of (\ref{eq:linear-sys}),
i.e.,
\begin{equation}
\boldsymbol{x}_{\text{PI}}=\arg\min_{\boldsymbol{x}}\left\Vert \boldsymbol{x}\right\Vert _{2}\qquad\text{subject to}\qquad\min\left\Vert \boldsymbol{b}-\boldsymbol{A}\boldsymbol{x}\right\Vert _{2},\label{eq:pi-sol}
\end{equation}
or equivalently, $\boldsymbol{x}_{\text{PI}}=\boldsymbol{A}^{+}\boldsymbol{b}$.
When $\boldsymbol{A}$ is large-scale, it is preferable to solve (\ref{eq:linear-sys})
using a \emph{Krylov subspace} (\emph{KSP}) method, such as GMRES
\cite{saad1986gmres}, which seeks a solution in the $k$th KSP 
\begin{equation}
\mathcal{K}_{k}\left(\boldsymbol{A},\boldsymbol{v}\right)=\text{span}\left\{ \boldsymbol{v},\boldsymbol{A}\boldsymbol{v},\ldots,\boldsymbol{A}^{k-1}\boldsymbol{v}\right\} \label{eq:ksp-def}
\end{equation}
at the $k$th iteration, where $\boldsymbol{v}$ is typically equal
to $\boldsymbol{b}$. It is well known that KSP methods can benefit
from robust and effective preconditioners for ill-conditioned problems.
This work introduces a software package called \emph{HIFIR}, which
delivers robust and computationally efficient preconditioners for
singular systems. As a side product, HIFIR also improves the robustness
in preconditioning ill-conditioned systems.

Compared to nonsingular systems, preconditioning (nearly) singular
systems is a very challenging problem. Several software packages offer
fairly robust and easy-to-use preconditioners for nonsingular systems,
such as the multilevel ILU (MLILU) in ILUPACK \cite{bollhofer2011ilupack},
the supernodal ILU in SuperLU \cite{li2011supernodal}, and various
parallel preconditioners in PETSc \cite{petsc-user-ref}. Conceptually,
such software packages construct a preconditioner $\boldsymbol{M}\in\mathbb{C}^{n\times n}$
that approximates $\boldsymbol{A}$ in that $\boldsymbol{M}^{-1}\approx\boldsymbol{A}^{-1}$.
Given $\boldsymbol{M}$, a right-preconditioned\footnote{We consider only right preconditioning because left preconditioning
alters the computation of residual vector $\boldsymbol{r}=\boldsymbol{b}-\boldsymbol{A}\boldsymbol{x}$
and in turn may lead to false stagnation or early terminations for
ill-conditioned systems \cite{ghai2019comparison}.} KSP method seeks a solution to 
\begin{equation}
\boldsymbol{A}\boldsymbol{M}^{-1}\boldsymbol{y}=\boldsymbol{b}\label{eq:right-preconditioned}
\end{equation}
in $\mathcal{K}_{k}\left(\boldsymbol{A}\boldsymbol{M}^{-1},\boldsymbol{v}\right)$,
where typically $\boldsymbol{v}=\boldsymbol{b}$, and then $\boldsymbol{x}=\boldsymbol{M}^{-1}\boldsymbol{y}$.
Ideally, the preconditioned KSP methods would converge significantly
faster than the unpreconditioned ones. However, when $\boldsymbol{A}$
is (nearly) singular and the system (\ref{eq:linear-sys}) is inconsistent,
there is a lack of robust algorithms and software. Some earlier techniques
used a CGLS-type KSP method (e.g., \cite{bjorck1996numerical,fong2011lsmr,paige1982lsqr}),
which is mathematically equivalent to solving the normal equation
using CG \cite{bjorck1996numerical,hestenes1952methods} or MINRES
\cite{paige1975solution,fong2011lsmr}. Those KSP methods tend to
converge slowly due to the squaring of the condition number by the
normal equation in the corresponding KSP \cite{jiao2021approximate}.
More recently, there has been significant interest in preconditioning
GMRES-type methods for singular systems or least-squares problems
\cite{hayami2010gmres,jiao2021approximate,morikuni2015convergence}.
For example, the so-called AB-GMRES \cite{hayami2010gmres} solves
the system $\boldsymbol{A}\boldsymbol{B}\boldsymbol{y}=\boldsymbol{b}$
using GMRES with $\mathcal{K}_{k}\left(\boldsymbol{A}\boldsymbol{B},\boldsymbol{b}\right)$,
and then $\boldsymbol{x}=\boldsymbol{B}\boldsymbol{y}$. Here, $\boldsymbol{B}$
plays a similar role as $\boldsymbol{M}^{-1}$, except that $\boldsymbol{B}$
may be singular (or rank deficient if $\boldsymbol{A}$ is rectangular).
Hayami et al. \cite{hayami2010gmres} constructed $\boldsymbol{B}$
based on robust incomplete factorization (RIF) of Benzi and T\r{u}ma
\cite{benzi2003robustpd,benzi2003robust}, which was originally developed
for CGLS-type methods. Although RIF could accelerate the convergence
of AB-GMRES in \cite{hayami2010gmres}, it was not robust in general
\cite{morikuni2015convergence}. In more recent works \cite{gould2017state,morikuni2015convergence},
$\boldsymbol{B}$ in AB-GMRES is typically chosen to be $\boldsymbol{A}^{H}$
(or $\boldsymbol{A}^{T}$ for real matrices), which unfortunately
squares the condition number (analogous to CGLS) and in turn can slow
down the convergence. This work aims to deliver a right preconditioner
that is more efficient and robust than RIF, and more importantly,
enables near-optimal convergence rates. We achieve this goal by leveraging
the new theory of \emph{$\epsilon$-accurate} \emph{approximated generalized
inverse} (\emph{AGI}) \cite{jiao2021approximate}, as outlined in
Section~\ref{sec:theoretical-foundation}.

Our development of HIFIR was based on our earlier software package
called \emph{HILUCSI} \cite{chen2021hilucsi}, which was a prototype
implementation of an MLILU for nonsingular saddle-point problems.
Compared to single-level ILUs (such as ILU($k$) and ILUTP, etc. \cite[Chapter 10]{saad2003iterative}),
MLILU is generally more robust for nonsingular indefinite systems
\cite{ghai2019comparison,chen2021hilucsi}. HILUCSI leveraged several
techniques in a novel way to achieve superior efficiency and robustness
for saddle-point problems than other MLILU libraries (such as ARMS
\cite{saad2002arms}, ILUPACK \cite{bollhofer2011ilupack}, and ILU++
\cite{mayer2007ilu++}). In particular, HILUCSI achieved high efficiency
by introducing a scalability-oriented dropping in a dual-thresholding
strategy in the fan-in ILU\footnote{The technique is also known as the Crout version of ILU \cite{li2003crout}
or left-looking \cite{eisenstat1981algorithms}, but we adopt the
terminology of ``fan in'' update, which is commonly used in parallel
computing \cite{demmel1993parallel} and is more suggestive in terms
of its algorithmic behavior.} for linear-time complexity in its factorization and solve. It improved
robustness by leveraging mixed symmetric and unsymmetric preprocessing
techniques at different levels and combining static and dynamic permutations.
HIFIR inherits some of these core algorithmic components of HILUCSI,
as described in Section~\ref{sec:Algorithmic-components}. However,
as an MLILU technique, HILUCSI was not robust for singular systems,
for example, when its final level is singular. HIFIR is designed to
overcome this deficiency by leveraging a rank-revealing factorization
in its final level, introducing iterative refinement to build a variable
preconditioner, and introducing a new pivoting strategy in its ILU
portion, as we will detail in Section~\ref{sec:Algorithmic-components}.

The main contributions of this work are as follows. First and foremost,
we introduce one of the first software libraries to improve the robustness
for ill-conditioned and (nearly) singular systems to achieve near
machine precision. Our software library, called HIFIR, or \emph{Hybrid
Incomplete Factorization with Iterative Refinement}, computes an AGI
\cite{jiao2021approximate} by hybridizing incomplete LU and rank-revealing
QR (RRQR) in a multilevel fashion. When used as a right-preconditioner
for GMRES, this hybridization enables (near) optimal convergence for
consistent or ill-conditioned systems. When fortified with iterative
refinement in FGMRES \cite{saad1993flexible}, HIFIR enables the robust
computation of the left null space and the pseudoinverse solution
of inconsistent systems. We have implemented HIFIR using template-based
objective-oriented programming in C++. For user-friendliness, the
C++ HIFIR library is header-only, with easy-to-use high-level interfaces
in MATLAB and Python. The software is open-source and has been made
available at \url{https://github.com/hifirworks/hifir}. Second, this
work also introduces a novel inverse-based rook pivoting (IBRP) in
the fan-in ILU. We describe efficient data structures for the efficient
implementation of IBRP and show that this new pivoting strategy improves
the robustness and efficiency for some challenging singular systems.
Third, HIFIR offers some advanced features, such as the support of
complex arithmetic, the ability to precondition both $\boldsymbol{A}$
and $\boldsymbol{A}^{H}$ using the same factorization, and the ability
to multiply by an AGI of the preconditioning operator. These features
enable the use of HIFIR as building blocks for advanced preconditioners,
such as (parallel) block preconditioners. In addition, HIFIR supports
mixed precision (e.g., double precision combined with single or potentially
half precision) for the input matrix and the preconditioner, which
is beneficial for heterogeneous hardware platforms and limited-memory
settings.

The remainder of this paper is organized as follows. In Section~\ref{sec:theoretical-foundation},
we give an overview of the theoretical foundation of HIFIR, including
optimality conditions, treatment for singular systems, etc. In Section~\ref{sec:Algorithmic-components},
we describe the algorithmic components of HIF and highlight some implementation
details. Section~\ref{sec:Applying-HIFIR} describes how to apply
HIF as a preconditioner, including iterative refinement. In Section~\ref{sec:software-library},
we introduce the application programming interfaces of HIFIR in C++,
MATLAB, and Python with example implementations. Section~\ref{sec:Illustration-of-HIFIR}
demonstrates HIFIR for some large-scale applications with indefinite
ill-conditioned and singular inconsistent systems. Finally, Section~\ref{sec:Conclusion}
concludes the paper with a discussion on future directions. For completeness,
we present the details of our data structures in Appendix~\ref{sec:Flexible-array-based}
and the complexity analysis of IBRP in Appendix~\ref{sec:Time-complexity-rook-pivoting}.

\section{\label{sec:theoretical-foundation}Theoretical foundations}

In this section, we give an overview of the theoretical foundations
of HIFIR. Most of the theory was based on that in \cite{jiao2021approximate},
except that we generalize the results from real matrices to complex
ones. We present some of the most relevant theoretical results for
completeness, but we omit the proofs because they follow the same
arguments as those in \cite{jiao2021approximate}.

\subsection{\label{subsec:optimal-rpo}Mathematically optimal right-preconditioning
operators for consistent systems}

Let us first consider the issue of optimal preconditioning for a consistent
system (\ref{eq:linear-sys}), where $\boldsymbol{b}$ is in the range
of $\boldsymbol{A}$ (i.e., $\boldsymbol{b}\in\mathcal{R}(\boldsymbol{A})$).
In floating-point arithmetic, the convergence rate of a KSP method
for such systems depends on the following generalized notion of condition
numbers.
\begin{defn}
\label{def:condition-number} Given a potentially singular matrix
$\boldsymbol{A}\in\mathbb{C}^{m\times n}$, the \emph{2-norm condition
number }of $\boldsymbol{A}$ is the ratio between the largest and
the smallest nonzero singular values of $\boldsymbol{A}$, i.e., $\kappa(\boldsymbol{A})=\sigma_{1}(\boldsymbol{A})/\sigma_{r}(\boldsymbol{A})$,
where $r=\text{rank}(\boldsymbol{A})$.
\end{defn}
To accelerate a KSP method for such a system, we solve a right-preconditioned
system
\begin{equation}
\boldsymbol{A}\boldsymbol{G}\boldsymbol{y}=\boldsymbol{b},\label{eq:prec-sys}
\end{equation}
using a KSP method, and then $\boldsymbol{x}=\boldsymbol{G}\boldsymbol{y}$.
We refer to $\boldsymbol{G}$ as a \emph{right preconditioning operator
}(\emph{RPO}). Ideally, we would like $\kappa(\boldsymbol{A}\boldsymbol{G})\ll\kappa(\boldsymbol{A})$.
For nonsingular systems, $\boldsymbol{G}$ is equivalent to $\boldsymbol{M}^{-1}$
in (\ref{eq:right-preconditioned}); for singular systems, $\boldsymbol{G}$
generalizes $\boldsymbol{M}^{-1}$. The symbol $\boldsymbol{G}$ signifies
that it is based on a \emph{generalized inverse}.
\begin{defn}
\label{def:gen-inverse}\cite[Definitions 2.2]{rao1972generalized}
Given a potentially rank-deficient $\boldsymbol{A}\in\mathbb{C}^{m\times n}$,
$\boldsymbol{A}^{g}$ is a \emph{generalized inverse} of $\boldsymbol{A}$
if and only if $\boldsymbol{A}\boldsymbol{A}^{g}\boldsymbol{A}=\boldsymbol{A}$.
\end{defn}
It is worth noting that the Moore--Penrose pseudoinverse $\boldsymbol{A}^{+}$
is a special case of generalized inverses. Although it might be tempting
to construct the RPO $\boldsymbol{G}$ to approximate $\boldsymbol{A}^{+}$,
the pseudoinverse is overly restrictive. The following two properties
of generalized inverses make them particularly relevant to right-preconditioning
singular systems.
\begin{prop}
\label{prop:decomp-ag}\cite[Proposition 3.4]{jiao2021approximate}
If $\boldsymbol{A}^{g}$ is a generalized inverse of $\boldsymbol{A}\in\mathbb{C}^{m\times n}$,
then $\boldsymbol{A}\boldsymbol{A}^{g}$ is diagonalizable, and its
eigenvalues are all zeros and ones. In other words, there exists a
nonsingular matrix $\boldsymbol{X}\in\mathbb{C}^{m\times m}$, such
that $\boldsymbol{A}\boldsymbol{A}^{g}=\boldsymbol{X}\begin{bmatrix}\boldsymbol{I}_{r}\\
 & \boldsymbol{0}
\end{bmatrix}\boldsymbol{X}^{-1}$, where $\boldsymbol{I}_{r}$ is the $r\times r$ identity matrix
with $r=\text{rank}(\boldsymbol{A})$. Conversely, if there exists
a nonsingular $\boldsymbol{X}$ such that $\boldsymbol{A}\boldsymbol{G}=\boldsymbol{X}\begin{bmatrix}\boldsymbol{I}_{r}\\
 & \boldsymbol{0}
\end{bmatrix}\boldsymbol{X}^{-1}$ for $r=\text{rank}(\boldsymbol{A})$, then $\boldsymbol{G}$ is a
generalized inverse of $\boldsymbol{A}$.
\end{prop}

\begin{prop}
\label{prop:condition-number}\cite[Proposition A.1]{jiao2021approximate}
Given $\boldsymbol{A}\in\mathbb{C}^{n\times n}$ of rank $r=\text{rank}(\boldsymbol{A})$
and a generalized inverse $\boldsymbol{A}^{g}$ with $\boldsymbol{A}\boldsymbol{A}^{g}=\boldsymbol{X}\begin{bmatrix}\boldsymbol{I}_{r}\\
 & \boldsymbol{0}
\end{bmatrix}\boldsymbol{X}^{-1}$, the condition number of $\boldsymbol{A}\boldsymbol{A}^{g}$ is bounded
by $\kappa(\boldsymbol{X})$, i.e., $\kappa(\boldsymbol{A}\boldsymbol{A}^{g})\leq\kappa(\boldsymbol{X})$.
\end{prop}
From Proposition~\ref{prop:decomp-ag}, it is easy to show that any
generalized inverse $\boldsymbol{A}^{g}$ (or a nonzero scalar multiple
of $\boldsymbol{A}^{g}$) enables a mathematically optimal RPO for
consistent systems in the following sense.
\begin{thm}
\label{thm:opt-consistent}\cite[Theorem 3.6]{jiao2021approximate}
Given $\boldsymbol{A}\in\mathbb{C}^{n\times n}$ and a generalized
inverse $\boldsymbol{A}^{g}$, then GMRES with RPO $\alpha\boldsymbol{A}^{g}$
with $\alpha\ne0$ converges to a least-squares solution $\boldsymbol{x}_{\text{LS}}$
of (\ref{eq:linear-sys}) after one iteration for all $\boldsymbol{b}\in\mathcal{R}(\boldsymbol{A})$
and $\boldsymbol{x}_{0}\in\mathbb{C}^{n}$. Conversely, if GMRES with
RPO $\boldsymbol{G}$ converges to a least-squares solution $\boldsymbol{x}_{\text{LS}}$
in one iteration for all $\boldsymbol{b}\in\mathcal{R}(\boldsymbol{A})$
and $\boldsymbol{x}_{0}\in\mathbb{C}^{n}$, then $\boldsymbol{G}$
is a scalar multiple of a generalized inverse of $\boldsymbol{A}$.
\end{thm}
A corollary of Theorem~\ref{thm:opt-consistent} is that if $\mathcal{R}(\boldsymbol{A}^{g})=\mathcal{R}(\boldsymbol{A}^{H})$,
then the computed $\boldsymbol{x}_{\text{LS}}$ is the pseudoinverse
solution of (\ref{eq:linear-sys}). Let $\Pi_{\mathcal{R}(\boldsymbol{A}^{H})}$
denote a projection onto $\boldsymbol{A}^{H}$. Given any $\boldsymbol{A}^{g}$,
it is easy to show that $\Pi_{\mathcal{R}(\boldsymbol{A}^{H})}\boldsymbol{A}^{g}$
is also a generalized inverse, and the computed $\boldsymbol{x}_{\text{LS}}$
with $\Pi_{\mathcal{R}(\boldsymbol{A}^{H})}\boldsymbol{A}^{g}$ as
the RPO is then the pseudoinverse solution. Theorem~\ref{thm:opt-consistent}
assumes exact arithmetic. With rounding errors, the condition number
of $\boldsymbol{A}\boldsymbol{A}^{g}$ must be bounded by a small
constant, which holds in general if $\kappa(\boldsymbol{X})$ is bounded
due to Proposition~\ref{prop:condition-number}.

Although the above theory may seem abstract, it suggests a new approach
for constructing RPO based on generalized inverses. In particular,
one option to construct $\boldsymbol{A}^{g}$ is to hybridize multilevel
ILU with a rank-revealing decomposition on its final Schur complement.
As an illustration, consider the following example that combines LU
factorization without pivoting with QR factorization with column pivoting
(QRCP) \cite{Golub13MC}.
\begin{example}
\label{exm:hybrid-factorization}Given $\boldsymbol{A}\in\mathbb{C}^{n\times n}$,
after $n_{1}$ steps of Gaussian elimination,
\begin{equation}
\boldsymbol{A}=\begin{bmatrix}\boldsymbol{L}_{1}\\
\boldsymbol{L}_{2} & \boldsymbol{I}_{n_{2}}
\end{bmatrix}\begin{bmatrix}\boldsymbol{I}_{n_{1}}\\
 & \boldsymbol{S}
\end{bmatrix}\begin{bmatrix}\boldsymbol{U}_{1} & \boldsymbol{U}_{2}\\
 & \boldsymbol{I}_{n_{2}}
\end{bmatrix},\label{eq:partialLU}
\end{equation}
where $\boldsymbol{S}\in\mathbb{C}^{n_{2}\times n_{2}}$ is the \emph{Schur
complement}. Clearly, $\text{rank}(\boldsymbol{A})=n_{1}+\text{rank}(\boldsymbol{S})$.
Let the QRCP of $\boldsymbol{S}$ be 
\begin{equation}
\boldsymbol{S}\boldsymbol{P}=\boldsymbol{Q}\begin{bmatrix}\boldsymbol{R}_{1} & \boldsymbol{R}_{2}\\
 & \boldsymbol{0}
\end{bmatrix},\label{eq:QRCP-Schur-Complement}
\end{equation}
where $\boldsymbol{Q}\in\mathbb{C}^{n_{2}\times n_{2}}$ is unitary
and $\boldsymbol{R}_{1}\in\mathbb{C}^{s\times s}$ for $s=\text{rank}(\boldsymbol{S})$
(i.e., the (numerical) rank of $\boldsymbol{S}$). Let $\hat{\boldsymbol{P}}$
and $\hat{\boldsymbol{Q}}$ be composed of the first $s$ columns
of $\boldsymbol{P}$ and $\boldsymbol{Q}$, respectively. Then, $\boldsymbol{S}^{g}=\hat{\boldsymbol{P}}\boldsymbol{R}_{1}^{-1}\hat{\boldsymbol{Q}}^{H}$
is a generalized inverse of $\boldsymbol{S}$ with $\boldsymbol{S}\boldsymbol{S}^{H}=\hat{\boldsymbol{Q}}\hat{\boldsymbol{Q}}^{H}=\boldsymbol{Q}\boldsymbol{I}_{s}\boldsymbol{Q}^{H}$.
Furthermore, 
\begin{equation}
\boldsymbol{G}=\begin{bmatrix}\boldsymbol{U}_{1} & \boldsymbol{U}_{2}\\
 & \boldsymbol{I}_{n_{2}}
\end{bmatrix}^{-1}\begin{bmatrix}\boldsymbol{I}_{n_{1}}\\
 & \boldsymbol{S}^{g}
\end{bmatrix}\begin{bmatrix}\boldsymbol{L}_{1}\\
\boldsymbol{L}_{2} & \boldsymbol{I}_{n_{2}}
\end{bmatrix}^{-1}=\begin{bmatrix}\boldsymbol{U}_{1} & \boldsymbol{U}_{2}\\
 & \boldsymbol{I}_{n_{2}}
\end{bmatrix}^{-1}\begin{bmatrix}\boldsymbol{I}_{n_{1}}\\
 & \boldsymbol{P}\begin{bmatrix}\boldsymbol{R}_{1}^{-1} & \boldsymbol{0}\\
\boldsymbol{} & \boldsymbol{0}
\end{bmatrix}\boldsymbol{Q}^{H}
\end{bmatrix}\begin{bmatrix}\boldsymbol{L}_{1}\\
\boldsymbol{L}_{2} & \boldsymbol{I}_{n_{2}}
\end{bmatrix}^{-1}\label{eq:mlilu-pseudoinverse}
\end{equation}
is a generalized inverse of $\boldsymbol{A}$, with $\boldsymbol{X}=\begin{bmatrix}\boldsymbol{L}_{1}\\
\boldsymbol{L}_{2} & \boldsymbol{I}_{n_{2}}
\end{bmatrix}\begin{bmatrix}\boldsymbol{I}_{n_{1}}\\
 & \boldsymbol{Q}
\end{bmatrix}$ for $\boldsymbol{X}$ as in Proposition~\ref{prop:decomp-ag}.
\end{example}
We refer to the preceding construction of $\boldsymbol{G}$ as a \emph{hybrid
factorization}. It enables a more efficient approach to construct
an optimal RPO, for example, compared to applying QRCP to $\boldsymbol{A}$
if $n_{2}\ll n$. For the strategy to be successful, we must address
some practical issues. First, we need to allow droppings in the factorization
to reduce computational cost and memory requirement, especially for
larger-scale systems. Second, we need to control $\kappa\left(\begin{bmatrix}\boldsymbol{L}_{1}\\
\boldsymbol{L}_{2} & \boldsymbol{I}
\end{bmatrix}\right)$ to limit $\kappa(\boldsymbol{A}\boldsymbol{G})$, for example, by
leveraging pivoting and equilibration \cite{duff2001algorithms}.
Third, it is desirable to make $\boldsymbol{S}$ as small as possible
before we apply QRCP. Hereafter, we will address the first two issues
from a theoretical perspective and then address the third issue in
Section~\ref{sec:Algorithmic-components}.

\subsection{Near-optimal right-preconditioning operators via approximate generalized
inverses}

Although a mathematically optimal RPO enables the most rapid convergence
of KSP methods, the computational cost per iteration may be prohibitively
high, so is the memory requirement. In practice, it may be more desirable
to construct ``near-optimal'' RPOs by approximating a generalized
inverse. The following definition and theorem establish the guideline
for constructing such approximations.
\begin{defn}
\label{def:epsilon-accuracy}\cite[Definition 3.8]{jiao2021approximate}
Given $\boldsymbol{A}\in\mathbb{C}^{n\times n}$, $\boldsymbol{G}$
is an \emph{$\epsilon$-accurate approximate generalized inverse}
(\emph{AGI}) if there exists $\boldsymbol{X}\in\mathbb{C}^{n\times n}$
such that
\begin{equation}
\left\Vert \boldsymbol{X}^{-1}\boldsymbol{A}\boldsymbol{G}\boldsymbol{X}-\begin{bmatrix}\boldsymbol{I}_{r}\\
 & \boldsymbol{0}
\end{bmatrix}\right\Vert =\epsilon\le1,\label{eq:epsilon-accuracy}
\end{equation}
where $\boldsymbol{I}_{r}$ is $r\times r$ identity matrix with $r=\text{rank}(\boldsymbol{A})$.
A class of AGI is $\epsilon$-accurate if $\epsilon$ tends to 0 as
its control parameters are tightened. $\boldsymbol{G}$ is a stable
AGI if $\kappa(\boldsymbol{X})\le C$ for some $C\ll\nicefrac{1}{\epsilon_{\text{mach}}}$.
\end{defn}
\begin{thm}
\label{thm:agi-conv-consistent}\cite[Theorem 3.9]{jiao2021approximate}
GMRES with an $\epsilon$-accurate AGI $\boldsymbol{G}$ of $\boldsymbol{A}$
converges to a least-squares solution $\boldsymbol{x}_{\text{LS}}$
of (\ref{eq:linear-sys}) in exact arithmetic for all consistent systems
(i.e., $\boldsymbol{b}\in\mathcal{R}(\boldsymbol{A})$) with any initial
guess $\boldsymbol{x}_{0}\in\mathbb{C}^{n}$.
\end{thm}
A corollary of Theorem~\ref{thm:agi-conv-consistent} is that given
an $\epsilon$-accurate AGI $\boldsymbol{G}$, $\Pi_{\mathcal{R}(\boldsymbol{A}^{H})}\boldsymbol{G}$
is also an $\epsilon$-accurate AGI, and the computed $\boldsymbol{x}_{\text{LS}}$
with $\Pi_{\mathcal{R}(\boldsymbol{A}^{H})}\boldsymbol{G}$ as the
RPO is the pseudoinverse solution. Mathematically, it is equivalent
to compute $\boldsymbol{x}_{\text{LS}}$ with $\boldsymbol{G}$ as
the RPO and then project $\boldsymbol{x}_{\text{LS}}$ onto $\mathcal{R}(\boldsymbol{A}^{H})$
to obtain the pseudoinverse solution.
\begin{rem}
In the literature, commonly used measures of accuracy and stability
of a preconditioner $\boldsymbol{M}$ for a nonsingular matrix $\boldsymbol{A}$
were $\Vert\boldsymbol{A}-\boldsymbol{M}\Vert_{F}$ and $\Vert\boldsymbol{I}-\boldsymbol{A}\boldsymbol{M}^{-1}\Vert_{F}$,
respectively; see, e.g., \cite{benzi2002preconditioning}. Our new
definitions of accuracy and stability in Definition~\ref{def:epsilon-accuracy}
are more general in that they apply to singular systems. In addition,
they are more rigorous in that they are based on Theorem~\ref{thm:agi-conv-consistent}
and Proposition~\ref{prop:condition-number}, respectively, instead
of based on empirical evidence \cite{benzi2002preconditioning}.
\end{rem}
Theorem~\ref{thm:agi-conv-consistent} and Example~\ref{exm:hybrid-factorization}
suggest that we can construct AGIs by replacing the LU factorization
in hybrid factorization with some ILU variants. We will refer to the
combination of ILU with a rank-revealing factorization on the Schur
complement as a \emph{hybrid incomplete factorization} (\emph{HIF}).
From the perspective of AGI, a good candidate ILU should satisfy three
critical criteria. First, the ILU needs to have prudent dropping strategies
to make the approximation as accurate as possible. Second, we must
be able to control $\kappa(\boldsymbol{A}\boldsymbol{G})$ effectively
for stability. Third, the computational cost and storage requirement
should ideally scale linearly (or near-linearly) with respect to the
input size. In the ILU literature \cite{chow1997experimental,saad2003iterative,bollhofer2006multilevel},
there had been significant attention to the first criterion. However,
the second criterion excludes simple ILU techniques without pivoting,
such as ILU($k$) \cite{saad2003iterative}. The third criterion excludes
ILU with relatively simple pivoting strategies, such as ILUTP \cite{chow1997experimental,saad2003iterative}
and its supernodal variants \cite{li2011supernodal}, which suffer
from superlinear complexity \cite{ghai2019comparison,chen2021hilucsi}.

Although the second and third criteria may appear self-contradicting
for traditional ILU techniques, they can be met by a well-designed
MLILU technique. Before delving into the details of MLILU algorithms,
let us briefly review MLILU and more importantly, show that MLILU
can be used to construct accurate and stable AGIs. First, consider
a two-level ILU (or more precisely, ILDU) of $\boldsymbol{A}\in\mathbb{C}^{n\times n}$,
\begin{equation}
\boldsymbol{P}^{T}\boldsymbol{W}\boldsymbol{A}\boldsymbol{V}\boldsymbol{Q}=\begin{bmatrix}\boldsymbol{B} & \boldsymbol{F}\\
\boldsymbol{E} & \boldsymbol{C}
\end{bmatrix}\approx\tilde{\boldsymbol{M}}=\begin{bmatrix}\tilde{\boldsymbol{B}} & \tilde{\boldsymbol{F}}\\
\tilde{\boldsymbol{E}} & \boldsymbol{C}
\end{bmatrix}=\begin{bmatrix}\boldsymbol{L}_{B}\\
\boldsymbol{L}_{E} & \boldsymbol{I}
\end{bmatrix}\begin{bmatrix}\boldsymbol{D}_{B}\\
 & \boldsymbol{S}
\end{bmatrix}\begin{bmatrix}\boldsymbol{U}_{B} & \boldsymbol{U}_{F}\\
 & \boldsymbol{I}
\end{bmatrix},\label{eq:two-lvl-ilu}
\end{equation}
where $\boldsymbol{B}\approx\tilde{\boldsymbol{B}}=\boldsymbol{L}_{B}\boldsymbol{D}_{B}\boldsymbol{U}_{B}$
is an ILDU of the leading block, $\boldsymbol{E}\approx\tilde{\boldsymbol{E}}=\boldsymbol{L}_{E}\boldsymbol{D}_{B}\boldsymbol{U}_{B}$,
$\boldsymbol{F}\approx\tilde{\boldsymbol{F}}=\boldsymbol{L}_{B}\boldsymbol{D}_{B}\boldsymbol{U}_{F}$,
and 
\begin{equation}
\boldsymbol{S}=\boldsymbol{C}-\tilde{\boldsymbol{E}}\tilde{\boldsymbol{B}}^{-1}\tilde{\boldsymbol{F}}=\boldsymbol{C}-\boldsymbol{L}_{E}\boldsymbol{D}_{B}\boldsymbol{U}_{F}\label{eq:schur-complement}
\end{equation}
is the Schur complement; $\boldsymbol{P}$ and $\boldsymbol{Q}$ are
row and column permutation matrices, respectively; $\boldsymbol{W}$
and $\boldsymbol{V}$ correspond to row and column scaling diagonal
matrices, respectively. The Schur complement $\boldsymbol{S}$ can
be factorized recursively using the same ILU technique, leading to
an $m$-level ILU preconditioners, namely,
\begin{equation}
\boldsymbol{M}=\underset{\boldsymbol{L}}{\underbrace{\boldsymbol{W}_{1}^{-1}\boldsymbol{P}_{1}\boldsymbol{L}_{1}\cdots\boldsymbol{W}_{m}^{-1}\boldsymbol{P}_{m}\boldsymbol{L}_{m}}}\begin{bmatrix}\boldsymbol{D}\\
 & \boldsymbol{S}_{m}
\end{bmatrix}\underset{\boldsymbol{U}}{\underbrace{\boldsymbol{U}_{m}\boldsymbol{Q}_{m}^{T}\boldsymbol{V}_{m}^{-1}\cdots\boldsymbol{U}_{1}\boldsymbol{Q}_{1}^{T}\boldsymbol{V}_{1}^{-1}}},\label{eq:multilevel-ilu}
\end{equation}
where $\boldsymbol{L}_{i}=\begin{bmatrix}\boldsymbol{I}_{n-n_{i}}\\
 & \boldsymbol{L}_{B}^{(i)}\\
 & \boldsymbol{L}_{E}^{(i)} & \boldsymbol{I}
\end{bmatrix}\in\mathbb{C}^{n\times n}$ for $\begin{bmatrix}\boldsymbol{L}_{B}^{(i)}\\
\boldsymbol{L}_{E}^{(i)} & \boldsymbol{I}
\end{bmatrix}\in\mathbb{C}^{n_{i}\times n_{i}}$ in (\ref{eq:two-lvl-ilu}) at the $i$th level (similarly for $\boldsymbol{U}_{i}$
and all other permutation and scaling matrices), $\boldsymbol{D}$
is composed of the ``union'' of $\boldsymbol{D}_{B}$ in (\ref{eq:two-lvl-ilu})
for all levels, and $\boldsymbol{S}_{m}$ is the final Schur complement.
As in Example~\ref{exm:hybrid-factorization}, we apply QRCP to $\boldsymbol{S}_{m}$
to obtain $\boldsymbol{S}_{m}\boldsymbol{P}=\boldsymbol{Q}\begin{bmatrix}\boldsymbol{R}_{1} & \boldsymbol{R}_{2}\\
 & \boldsymbol{0}
\end{bmatrix}$, where $\boldsymbol{Q}\in\mathbb{C}^{n_{m}\times n_{m}}$ is unitary
and $\boldsymbol{R}_{1}\in\mathbb{C}^{n_{s}\times n_{s}}$ with $n_{s}=\text{rank}(\boldsymbol{S}_{m})$,
$r_{11}\geq r_{22}\geq\dots\ge r_{n_{s}n_{s}}>0$ along its diagonal,
and $\kappa(\boldsymbol{R}_{1})\ll1/\epsilon_{\text{mach}}$. We define
an RPO as
\begin{equation}
\boldsymbol{G}=\boldsymbol{U}^{-1}\begin{bmatrix}\boldsymbol{D}^{-1}\\
 & \boldsymbol{S}_{m}^{g}
\end{bmatrix}\boldsymbol{L}^{-1}=\boldsymbol{U}^{-1}\begin{bmatrix}\boldsymbol{D}^{-1}\\
 & \boldsymbol{P}\begin{bmatrix}\boldsymbol{R}_{1}^{-1} & \boldsymbol{0}\\
\boldsymbol{} & \boldsymbol{0}
\end{bmatrix}\boldsymbol{Q}^{H}
\end{bmatrix}\boldsymbol{L}^{-1},\label{eq:hif-ig}
\end{equation}
where $\boldsymbol{L}$, $\boldsymbol{U}$, and $\boldsymbol{D}$
are as in (\ref{eq:multilevel-ilu}). We note the following fact.
\begin{prop}
\label{prop:mlilu-gi}\cite[Lemma 4.3]{jiao2021approximate} If no
dropping is applied in MLILU, then $\boldsymbol{G}$ in (\ref{eq:hif-ig})
is a generalized inverse of $\boldsymbol{A}$ with $\boldsymbol{X}=\boldsymbol{L}\begin{bmatrix}\boldsymbol{I}_{n-n_{m}}\\
 & \boldsymbol{Q}
\end{bmatrix}$ for $\boldsymbol{X}$ as in Proposition~\ref{prop:decomp-ag}.
\end{prop}
Since $\boldsymbol{G}$ without dropping is an optimal RPO for consistent
systems due to Theorem~\ref{thm:opt-consistent}, we claim that an
$\epsilon$-accurate and stable $\boldsymbol{G}$ constitutes a near-optimal
RPO for consistent systems. The near optimality requires sufficient
small droppings and numerical stability, or more precisely, $\boldsymbol{X}^{-1}\boldsymbol{A}\boldsymbol{G}\boldsymbol{X}$
should be close to $\begin{bmatrix}\boldsymbol{I}_{r}\\
 & \boldsymbol{0}
\end{bmatrix}$ for $r=\text{rank}(\boldsymbol{A})$ and $\kappa(\boldsymbol{L})$
must be controlled by the algorithm. To this end, we utilize an MLILU
technique called \emph{HILUCSI}, which stands for \emph{Hierarchical
Incomplete LU-Crout with Scalability-oriented and Inverse-based droppings}
\cite{chen2021hilucsi}. HILUCSI leverages several techniques, including
fan-in ILU \cite{li2003crout}, equilibration \cite{duff2001algorithms},
static and dynamic pivoting across different levels \cite{chen2021hilucsi,bollhofer2006multilevel},
etc., to meet the accuracy and stability requirements. As the name
suggests, HILUCSI focuses on scalability in terms of problem sizes,
and it has linear time complexity in each level for both the factorization
and solve stages \cite{chen2021hilucsi}. We defer the detailed description
of these algorithmic components to Section~\ref{sec:Algorithmic-components}. 
\begin{rem}
Besides HILUCSI, there were several MLILU software packages, such
as ARMS \cite{saad2002arms}, ILU++ \cite{mayer2006alternative},
ILUPACK \cite{bollhofer2011ilupack}, etc. Proposition~\ref{prop:mlilu-gi}
can also be applied to improve those packages to solve singular systems
by applying QRCP to the final Schur complement. To harness this benefit,
however, one must also extend them (especially ARMS and ILU++) to
ensure the stability of the $\boldsymbol{L}$ factor. 

Finally, we note that it is sometimes needed to reuse HIF to construct
preconditioners for both $\boldsymbol{A}$ and $\boldsymbol{A}^{H}$.
To achieve, we note the following property:
\end{rem}
\begin{prop}
\cite[Proposition 3.5]{jiao2021approximate} \label{prop:transpose-gen-inv}If
$\boldsymbol{A}^{g}$ is a generalized inverse of $\boldsymbol{A}\in\mathbb{C}^{m\times n}$,
then $\boldsymbol{A}^{gH}\equiv\left(\boldsymbol{A}^{g}\right)^{H}$
is a generalized inverse of $\boldsymbol{A}^{H}$.
\end{prop}
Conceptually, we can extend (\ref{eq:two-lvl-ilu}) to construct a
preconditioner $\tilde{\boldsymbol{M}}^{H}$ for $\boldsymbol{A}^{H}$
as 
\begin{equation}
\boldsymbol{Q}^{T}\boldsymbol{V}\boldsymbol{A}^{H}\boldsymbol{W}\boldsymbol{P}\approx\tilde{\boldsymbol{M}}^{H}=\begin{bmatrix}\tilde{\boldsymbol{B}}^{H} & \tilde{\boldsymbol{E}}^{H}\\
\tilde{\boldsymbol{F}}^{H} & \boldsymbol{C}^{H}
\end{bmatrix}=\begin{bmatrix}\boldsymbol{U}_{B}^{H}\\
\boldsymbol{U}_{F}^{H} & \boldsymbol{I}
\end{bmatrix}\begin{bmatrix}\bar{\boldsymbol{D}}_{B}\\
 & \boldsymbol{S}^{H}
\end{bmatrix}\begin{bmatrix}\boldsymbol{L}_{B}^{H} & \boldsymbol{L}_{E}^{H}\\
 & \boldsymbol{I}
\end{bmatrix}.\label{eq:conjugate-trans}
\end{equation}

\subsection{\label{subsec:Variable-preconditioning}Variable preconditioning
via iterative refinement for null-space computation}

The preceding discussions focused on consistent systems. For inconsistent
systems, i.e., $\boldsymbol{b}\in\mathbb{C}^{n}\backslash\mathcal{R}(\boldsymbol{A})$,
an AGI (or even a generalized inverse\footnote{When $\boldsymbol{G}=\boldsymbol{A}^{g}$, then right-preconditioned
GMRES converges to a weighted-least-squares (WLS) solution for inconsistent
systems, instead of least-squares solution \cite[Theorem 3.7]{jiao2021approximate}.
We cannot convert the WLS solution into a pseudoinverse solution by
projecting it onto $\mathcal{R}(\boldsymbol{A}^{H})$.}) cannot guarantee the convergence of GMRES due to the following fact.
\begin{thm}
\label{thm:range-symmetry}\cite[Theorem 2.4]{jiao2021approximate}
GMRES with RPO $\boldsymbol{G}$ does not break down until finding
a least-squares solution $\boldsymbol{x}_{\text{LS}}$ of (\ref{eq:linear-sys})
for all $\boldsymbol{b}\in\mathbb{C}^{n}$ and $\boldsymbol{x}_{0}\in\mathbb{C}^{n}$
if and only if $\boldsymbol{A}\boldsymbol{G}$ is range symmetric
(i.e., $\mathcal{R}(\boldsymbol{A}\boldsymbol{G})=\mathcal{R}((\boldsymbol{A}\boldsymbol{G})^{H})$)
and $\mathcal{R}(\boldsymbol{A})=\mathcal{R}(\boldsymbol{A}\boldsymbol{G})$.
Furthermore, $\boldsymbol{x}_{\text{LS}}$ is the pseudoinverse solution
if $\mathcal{R}(\boldsymbol{G})=\mathcal{R}(\boldsymbol{A}^{H})$.
\end{thm}
\begin{rem}
The main challenge posed by Theorem~\ref{thm:range-symmetry} is
that it is difficult, if not impossible, to ensure the range symmetry
of $\boldsymbol{A}\boldsymbol{G}$ for an approximate generalized
inverse $\boldsymbol{G}$. The requirement of range symmetry is the
primary reason why $\boldsymbol{A}^{H}$ is often used as $\boldsymbol{B}$
in AB-GMRES \cite{gould2017state,morikuni2015convergence}. It is
also a key factor for the prevalence of CGLS-type KSP methods \cite{bjorck1996numerical,fong2011lsmr,paige1982lsqr}
for solving singular and least-squares problems. Although such methods
can be accelerated by applying some preconditioners, such as incomplete
QR \cite{jennings1984incomplete,saad1988preconditioning} or RIF \cite{benzi2003robustpd,benzi2003robust},
it is difficult for these preconditioners to overcome the slowdown
caused by the squaring of the condition number by the normal equation.
\end{rem}
Fortunately, this issue can be resolved by using flexible GMRES with
\emph{variable preconditioners}, as formalized by the following definition
and theorem.
\begin{defn}
\label{def:fksp}\cite[Definition 2.1]{jiao2021approximate} Given
a matrix $\boldsymbol{A}\in\mathbb{C}^{n\times n}$, an initial vector
$\boldsymbol{v}\in\mathbb{C}^{n}$, and variable preconditioners $\boldsymbol{\mathcal{G}}_{k-1}=[\mathcal{G}_{1},\mathcal{G}_{2},\dots,\mathcal{G}_{k-1}]$,
the $k$th \emph{flexible Krylov subspace }(\emph{FKSP}) associated
with $\boldsymbol{A}$, $\boldsymbol{v}$, and $\boldsymbol{\mathcal{G}}_{k-1}$
is 
\begin{equation}
\mathcal{K}_{k}\left(\boldsymbol{A},\boldsymbol{v},\boldsymbol{\mathcal{G}}_{k-1}\right)=\text{span}\left\{ \boldsymbol{v},\boldsymbol{A}\mathcal{\mathcal{G}}_{1}(\boldsymbol{v}_{1}),\dots,\boldsymbol{A}\mathcal{\mathcal{G}}_{k-1}(\boldsymbol{v}_{k-1})\right\} ,\label{eq:FKSP}
\end{equation}
where $\boldsymbol{v}_{k-1}\in\mathcal{K}_{k-1}(\boldsymbol{A},\boldsymbol{v},\boldsymbol{\mathcal{G}}_{k-2})\backslash\{\boldsymbol{0}\}$
and $\boldsymbol{v}_{k-1}\perp\mathcal{K}_{k-2}(\boldsymbol{A},\boldsymbol{v},\boldsymbol{\mathcal{G}}_{k-3})$
for $k\geq2$. The \emph{flexible Krylov matrix, denoted by} $\boldsymbol{K}_{k}$,
is composed of the basis vectors in (\ref{eq:FKSP}).
\end{defn}
The orthogonality of $\boldsymbol{v}_{k-1}$ with $\mathcal{K}_{k-2}(\boldsymbol{A},\boldsymbol{v},\boldsymbol{\mathcal{G}}_{k-3})$
can be enforced using a generalization of Arnoldi iterations as in
\cite{saad1993flexible}.
\begin{thm}
\label{thm:inconsistent-cvg}\cite[Theorem 2.5]{jiao2021approximate}
If FGMRES with variable preconditioners $\boldsymbol{\mathcal{G}}=[\mathcal{G}_{1},\mathcal{G}_{2},\ldots]$
does not break down until step $k+1$, where $k=\text{rank}(\boldsymbol{A})$
for a specific $\boldsymbol{b}\in\mathcal{R}(\boldsymbol{A})$ and
$k=\text{rank}(\boldsymbol{A})+1$ for $\boldsymbol{b}\in\mathbb{C}^{n}\backslash\mathcal{R}(\boldsymbol{A})$,
respectively. If $\mathcal{R}(\boldsymbol{A})=\sum_{i=1}^{k}\mathcal{R}(\boldsymbol{A}\mathcal{G}_{i})$,
then it finds a least-squares solution of (\ref{eq:linear-sys}) for
an initial guess $\boldsymbol{x}_{0}\in\mathbb{C}^{n}$.
\end{thm}
One effective approach to construct variable preconditioners is to
introduce \emph{iterative refinement} (\emph{IR}) in HIF, leading
to HIFIR. Specifically, given an $\epsilon$-accurate AGI $\boldsymbol{G}$
and an initial vector $\boldsymbol{v}_{0}\in\mathbb{C}^{n}$ (typically,
$\boldsymbol{v}_{0}=\boldsymbol{0}$), we refine the solution of $\boldsymbol{A}\boldsymbol{v}=\boldsymbol{q}$
iteratively by obtaining $\boldsymbol{v}_{j}$ for $j=1,2,\ldots$
as
\begin{equation}
\boldsymbol{v}_{j}=\left(\boldsymbol{I}-\boldsymbol{G}\boldsymbol{A}\right)\boldsymbol{v}_{j-1}+\boldsymbol{G}\boldsymbol{q}.\label{eq:ir}
\end{equation}
 Eq.~(\ref{eq:ir}) defines the $i$th RPO $\mathcal{G}_{i}$ in
$\boldsymbol{\mathcal{G}}$, and different $\mathcal{G}_{i}$ may
use different numbers of IR iterations. Note that (\ref{eq:ir}) in
general does not converge by a standalone iterative solver in that
$\rho(\boldsymbol{I}-\boldsymbol{G}\boldsymbol{A})\ge1$ when $\boldsymbol{G}$
is a generalized inverse, and Eq.~(\ref{eq:ir}) cannot introduce
additional nonlinearity into the variable preconditioner $\boldsymbol{\mathcal{G}}$
and in turn undermine the robustness of FGMRES due to Theorem~\ref{thm:inconsistent-cvg}.
Hence, we shift our attention to its use as a variable preconditioner
in the context of computing null-space vectors. In particular, we
apply the technique to compute the left null space of $\boldsymbol{A}$
(i.e., $\mathcal{N}(\boldsymbol{A}^{H})$), which allows us to convert
an inconsistent system into a consistent one. Furthermore, by applying
the same technique to compute the right null space of $\boldsymbol{A}$
(i.e., $\mathcal{N}(\boldsymbol{A})$), we can convert a least-squares
solution from the consistent system into the pseudoinverse solution
\cite{jiao2021approximate}.

\section{\label{sec:Algorithmic-components}Recursive construction of hybrid
incomplete factorization}

\begin{figure}
\centering{}\includegraphics[width=1\columnwidth]{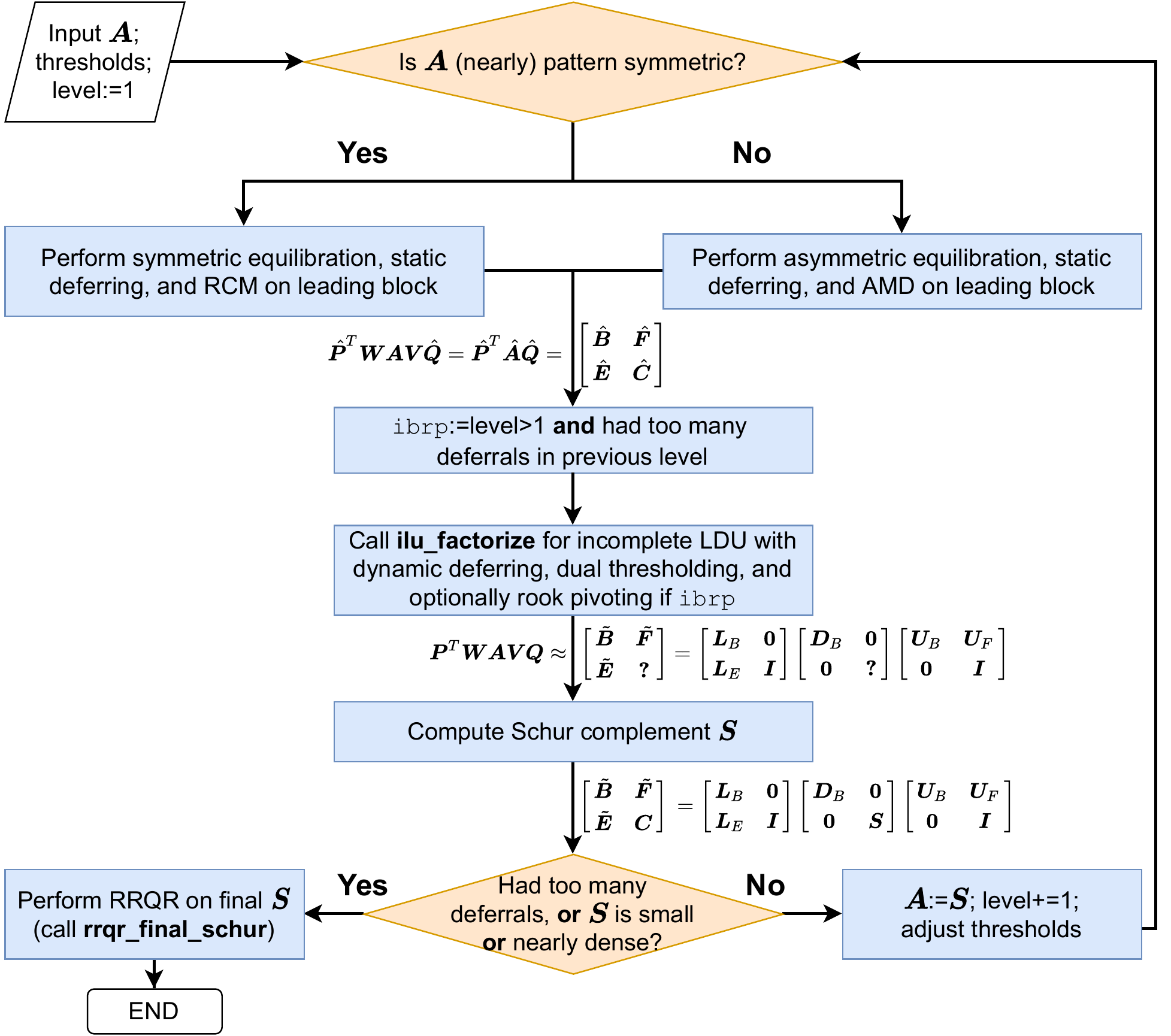}\caption{\label{fig:hif-workflow}Overall control flow of hybrid incomplete
factorization.}
\end{figure}

In this section, we describe the overall algorithm for constructing
HIF. Similar to that of HILUCSI \cite{chen2021hilucsi} and other
MLILU algorithms, HIF is a recursive algorithm. As shown in Figure~\ref{fig:hif-workflow},
at each level, HIF takes an input matrix $\boldsymbol{A}$. HIF first
performs symmetric or unsymmetric preprocessing depending on whether
$\boldsymbol{A}$ is (nearly) pattern symmetric, i.e., $\text{nnp}(\boldsymbol{A})\approx\text{nnp}(\boldsymbol{A}^{T})$;
see Section~\ref{subsec:Symmetric-and-asymmetric}. The process leads
to $\hat{\boldsymbol{P}}^{T}\boldsymbol{W}\boldsymbol{A}\boldsymbol{V}\hat{\boldsymbol{Q}}=\hat{\boldsymbol{P}}^{T}\hat{\boldsymbol{A}}\hat{\boldsymbol{Q}}$,
where $\hat{\boldsymbol{P}}$, $\boldsymbol{W}$, $\boldsymbol{V}$,
and $\hat{\boldsymbol{Q}}$ are obtained from preprocessing. If $\boldsymbol{A}$
is nearly pattern symmetric, the preprocessing step may involve static
deferring, which splits $\hat{\boldsymbol{P}}^{T}\hat{\boldsymbol{A}}\hat{\boldsymbol{Q}}$
into a 2-by-2 block matrix, i.e., $\hat{\boldsymbol{P}}^{T}\hat{\boldsymbol{A}}\hat{\boldsymbol{Q}}=\begin{bmatrix}\hat{\boldsymbol{B}} & \hat{\boldsymbol{F}}\\
\hat{\boldsymbol{E}} & \hat{\boldsymbol{C}}
\end{bmatrix}$. Depending on whether the leading block $\hat{\boldsymbol{B}}$ is
Hermitian, we then perform incomplete $\boldsymbol{L}\boldsymbol{D}\boldsymbol{L}^{H}$
or $\boldsymbol{L}\boldsymbol{D}\boldsymbol{U}$ factorizations, respectively,
where $\boldsymbol{L}$ and $\boldsymbol{U}$ are unit lower and upper
triangular matrices, respectively (i.e., their diagonal entries are
ones). We compute these factorizations using \emph{fan-in updates}
(aka the Crout version of ILU), which update the $k$th column of
$\boldsymbol{L}$ using columns 1 through $k-1$ and update the $k$th
row of $\boldsymbol{U}$ using rows 1 through $k-1$ at the $k$th
step, respectively. For stability, we combine fan-in updates with
dynamic deferring and scalability-oriented droppings; see Section~\ref{subsec:Incomplete-factorization}.
In addition, HIF enables rook pivoting when the previous level had
too many deferrals as indicated by the Boolean tag \textsf{ibrp} in
Figure~\ref{fig:hif-workflow}; see Section~\ref{subsec:Inverse-based-rook-pivoting}.
The dynamic deferring may permute some rows and columns in $\hat{\boldsymbol{B}}$
after $\hat{\boldsymbol{E}}$ and $\hat{\boldsymbol{F}}$ to obtain
a new 2-by-2 block structure $\begin{bmatrix}\tilde{\boldsymbol{B}} & \tilde{\boldsymbol{F}}\\
\tilde{\boldsymbol{E}} & \boldsymbol{C}
\end{bmatrix}$, where $\hat{\boldsymbol{E}}$ and $\hat{\boldsymbol{F}}$ are the
leading rows and blocks in $\boldsymbol{E}$ and $\boldsymbol{F}$,
respectively. We then compute the Schur complement $\boldsymbol{S}$
corresponding to $\boldsymbol{C}$ and factorize $\boldsymbol{S}$
either directly using RRQR or recursively using HIF, depending on
whether $\boldsymbol{S}$ is sufficiently small or nearly dense; see
Section~\ref{subsec:Schur-complement}. Some of the components above
are the same as those in HILUCSI \cite{chen2021hilucsi}, the predecessor
of HIF. The key differences between HIF and HILUCSI are that 1) HIF
introduces a variant of the rook pivoting \cite{poole2000rook} to
improve the stability of ILU and in turn, reduces the size of the
final Schur complement $\boldsymbol{S}_{m}$ and 2) HIF uses a rank-revealing
QR factorization \cite{chan1987rank} on $\boldsymbol{S}_{m}$. In
the following, we first focus on these two aspects and then briefly
outline the other components. We will describe how to apply HIF as
a preconditioner in Section~\ref{sec:Applying-HIFIR}.

\subsection{Incomplete $\boldsymbol{L}\boldsymbol{D}\boldsymbol{U}$ with dynamic
deferring and scalability-oriented and inverse-based droppings\label{subsec:Incomplete-factorization}}

The core of HIF at each level is incomplete $\boldsymbol{L}\boldsymbol{D}\boldsymbol{U}$
factorization, or its variant of incomplete $\boldsymbol{L}\boldsymbol{D}\boldsymbol{L}^{H}$
factorization. For robustness, we leverage the fan-in ILU, scalability-oriented
dropping, inverse-based dropping, and dynamic deferring, as shown
in Algorithm~\ref{alg:ilu_factor}. Hereafter, we outline the dual
thresholding and dynamic deferring steps in fan-in ILU. An optional
step in Algorithm~\ref{alg:ilu_factor} is the inverse-based rook
pivoting activated by the Boolean flag \textsf{ibrp}. We will dedicate
Section~\ref{subsec:Inverse-based-rook-pivoting} to this new pivoting
strategy.

In Algorithm~\ref{alg:ilu_factor}, the input matrix $\hat{\boldsymbol{A}}=\boldsymbol{W}\boldsymbol{A}\boldsymbol{V}$
is obtained from applying preprocessing techniques on either the input
matrix $\boldsymbol{A}$ or the Schur complement from the previous
level. Note that in the actual implementation, we do not form $\hat{\boldsymbol{A}}$
explicitly; instead, we rescale the entries in $\boldsymbol{A}$ by
$\boldsymbol{W}$ and $\boldsymbol{V}$ in a ``just-in-time'' fashion.
The procedure \textbf{ilu\_factorize} also takes $\hat{\boldsymbol{P}}$
and $\hat{\boldsymbol{Q}}$ are as input, which were obtained from
preprocessing along with $\boldsymbol{W}$ and $\boldsymbol{V}$ (see
Figure~\ref{fig:hif-workflow} and Section~\ref{subsec:Symmetric-and-asymmetric}).
The main loop in Algorithm~\ref{alg:ilu_factor} factorizes $\hat{\boldsymbol{B}}$,
the leading block of $\hat{\boldsymbol{A}}$, using fan-in ILU similar
to those in \cite{li2003crout,li2005crout}, which delays the computation
of the Schur complement as late as possible. Unlike \cite{li2003crout,li2005crout},
however, \textbf{ilu\_factorize} dynamically permutes (aka \emph{defers})
rows and columns in $\hat{\boldsymbol{B}}$ to the end of $\hat{\boldsymbol{E}}$
and $\hat{\boldsymbol{F}}$. The procedure computes 
\begin{equation}
\boldsymbol{P}^{T}\boldsymbol{W}\boldsymbol{A}\boldsymbol{V}\boldsymbol{Q}\approx\begin{bmatrix}\tilde{\boldsymbol{B}} & \tilde{\boldsymbol{F}}\\
\tilde{\boldsymbol{E}} & \boldsymbol{?}
\end{bmatrix}=\begin{bmatrix}\boldsymbol{L}_{B}\\
\boldsymbol{L}_{E} & \boldsymbol{I}
\end{bmatrix}\begin{bmatrix}\boldsymbol{D}_{B}\\
 & \boldsymbol{?}
\end{bmatrix}\begin{bmatrix}\boldsymbol{U}_{B} & \boldsymbol{U}_{F}\\
 & \boldsymbol{I}
\end{bmatrix},\label{eq:ilu_factor}
\end{equation}
where the first and second question marks (`$\boldsymbol{?}$') in
(\ref{eq:ilu_factor}) correspond to $\boldsymbol{C}$ and $\boldsymbol{S}$,
which we will describe their computations in Section~\ref{subsec:Schur-complement}.
Algorithm~\ref{alg:ilu_factor} returns $\vec{L}_{B}$, $\vec{D}_{B}$,
$\vec{U}_{B}$, $\vec{L}_{E}$, $\vec{U}_{F}$, $\vec{P}$, and $\vec{Q}$,
from which $\tilde{\boldsymbol{E}}$ and $\tilde{\boldsymbol{F}}$
(along with $\boldsymbol{C}$) can be computed.

A noteworthy feature of \textbf{ilu\_factorize} is its \emph{scalability-oriented
dropping}, which differs from the dropping strategies in \cite{li2003crout,li2005crout}
and \cite{bollhofer2011ilupack}. Consider the step $k$ of ILU at
a particular level. Using the MATLAB's colon notation, let $\boldsymbol{\ell}_{k+1:n,k}$
and $\boldsymbol{u}_{k,k+1:n}$ denote the $k$th column and row of
$\begin{bmatrix}\boldsymbol{L}_{B}\\
\boldsymbol{L}_{E}
\end{bmatrix}$ and $\begin{bmatrix}\boldsymbol{U}_{B} & \boldsymbol{U}_{F}\end{bmatrix}$,
respectively. Our dropping strategy limits $\text{nnz}\left(\boldsymbol{\ell}_{k+1:n,k}\right)$
and $\text{nnz}\left(\boldsymbol{u}_{k,k+1:n}\right)$ to be proportional
to numbers of nonzeros in the corresponding column and row of the
original (i.e., the top level instead of the present level) input
matrix, respectively; see line~\ref{line:ilu:select-ell-u}. This
dropping strategy plays an important role for HILUCSI and HIF to achieve
(near) linear complexity in both space and time. Besides this symbolic
dropping, HILUCSI and HIF also adopted an inverse-based dropping \cite{bollhofer2006multilevel},
which drops every entry $\ell_{ik}$ such that $\kappa_{D}\left\Vert \boldsymbol{L}_{1:k,1:k}^{-1}\right\Vert _{\infty}\left|\ell_{ik}\right|\le\tau_{L}$,
where $\kappa_{D}$ and $\tau_{L}$ are user-specified thresholds
for the upper bound of $\left\Vert \boldsymbol{D}_{1:k,1:k}^{-1}\right\Vert $
and the drop tolerance, respectively (see Section~\ref{subsec:Control-parameters}).
The dropping for $\boldsymbol{U}$ is similar; see line~\ref{line:ilu:num-drop}. 

\begin{algorithm}
\begin{raggedright}
\caption{\label{alg:ilu_factor}\textbf{ilu\_factorize}$\left(\hat{\boldsymbol{A}},\boldsymbol{p},\boldsymbol{q},n_{1},\text{params},\text{level},\text{ibrp},\text{\textbf{nr}},\text{\textbf{nc}}\right)$}
\textbf{inputs}:
\par\end{raggedright}
\begin{raggedright}
\hspace{1cm}\textbf{ }$\hat{\boldsymbol{A}}$: input scaled matrix
of size $n\times n$ (i.e., $\hat{\boldsymbol{A}}=\boldsymbol{W}\boldsymbol{A}\boldsymbol{V}$,
passed in as $\boldsymbol{A}$, $\boldsymbol{W}$, and $\boldsymbol{V}$
separately)
\par\end{raggedright}
\begin{raggedright}
\hspace{1cm} $\boldsymbol{p},\boldsymbol{q}$: row and column permutation
vectors of $\boldsymbol{A}$ after preprocessing, respectively
\par\end{raggedright}
\begin{raggedright}
\hspace{1cm} $n_{1}$: the dimension of the current leading block
(i.e., $\boldsymbol{B}=\hat{\boldsymbol{A}}_{\boldsymbol{p}_{1:n_{1}},\boldsymbol{q}_{1:n_{1}}}$)
\par\end{raggedright}
\begin{raggedright}
\hspace{1cm} $\text{params}$: $\alpha_{L}$, $\alpha_{U}$, $\kappa$,
$\kappa_{D}$, $\tau_{L}$, $\tau_{U}$ (adapted for present level),
max\_steps (for rook pivoting)
\par\end{raggedright}
\begin{raggedright}
\hspace{1cm} level: current level
\par\end{raggedright}
\begin{raggedright}
\hspace{1cm} ibrp: Boolean tag for enabling inverse-based rook pivoting
\par\end{raggedright}
\begin{raggedright}
\hspace{1cm} $\text{\textbf{nr}},\text{\textbf{nc}}$: number of
nonzeros per row and column entry of the original user input matrix,
respectively
\par\end{raggedright}
\begin{raggedright}
\textbf{outputs}:
\par\end{raggedright}
\begin{raggedright}
\hspace{1cm} $\boldsymbol{L}_{B},\boldsymbol{d}_{B},\boldsymbol{U}_{B}$:
approximate $\boldsymbol{L}\boldsymbol{D}\boldsymbol{U}$ factors
of the leading block
\par\end{raggedright}
\begin{raggedright}
\hspace{1cm} $\boldsymbol{L}_{E},\boldsymbol{U}_{F}$: off-diagonal
blocks of (\ref{eq:two-lvl-ilu})
\par\end{raggedright}
\begin{raggedright}
\hspace{1cm} $\boldsymbol{p},\boldsymbol{q}$: updated row and column
permutation vectors, respectively
\par\end{raggedright}
\begin{raggedright}
\hspace{1cm} $n_{1}$: updated leading block dimension
\par\end{raggedright}
\begin{algorithmic}[1]

\STATE $\boldsymbol{L}\leftarrow\begin{bmatrix}1\end{bmatrix}$;
$\boldsymbol{U}\leftarrow\begin{bmatrix}1\end{bmatrix}$; $\boldsymbol{d}\leftarrow[]$

\FOR{$k=1$ \TO $n_{1}$}

\IF{ibrp}

\STATE $\boldsymbol{p},\ \boldsymbol{L},\ \boldsymbol{d},\ \boldsymbol{U},\ \boldsymbol{q}\leftarrow\text{\textbf{ib\_rook\_pivot}}\left(\hat{\boldsymbol{A}},k,\boldsymbol{p},\boldsymbol{L},\boldsymbol{d},\boldsymbol{U},\boldsymbol{q},\kappa,n_{1},\text{max\_steps}\right)$\label{line:ilu:pivot}\hfill{}\{perform
inverse-based root pivoting, typically only if $\text{level}>1$\}

\ELSE 

\STATE $d_{k}\leftarrow a_{p_{k}q_{k}}-\boldsymbol{\ell}_{k,1:k-1}\boldsymbol{D}_{1:k-1,1:k-1}\boldsymbol{u}_{1:k-1,k}$\hfill{}\COMMENT{fan-in
update of diagonal entry}

\ENDIF

\STATE $\tilde{\kappa}_{L}\leftarrow\left\Vert \boldsymbol{L}_{1:k,1:k}^{-1}\right\Vert _{\infty}$;
$\tilde{\kappa}_{U}\leftarrow\left\Vert \boldsymbol{U}_{1:k,1:k}^{-1}\right\Vert _{1}$\label{line:ilu:inv-norm}\hfill{}\{estimate
inverse norms of current $\boldsymbol{L}$ and $\boldsymbol{U}$ factors\}

\WHILE{$\kappa_{D}\left|d_{k}\right|<1$ \OR $\max\left\{ \tilde{\kappa}_{L},\tilde{\kappa}_{U}\right\} >\kappa$}\label{line:ilu:poor-factor}

\STATE permute entries $\boldsymbol{\ell}_{k,1:k-1}$, $d_{k}$ and
$\boldsymbol{u}_{1:k-1,k}$ to the end; update $\boldsymbol{p}$ and
$\boldsymbol{q}$ accordingly\label{line:ilu:dyn-defer}

\STATE \textbf{break if }$k==n_{1}$

\STATE $\tilde{\kappa}_{L}\leftarrow\left\Vert \boldsymbol{L}_{1:k.1:k}^{-1}\right\Vert _{\infty}$;
$\tilde{\kappa}_{U}\leftarrow\left\Vert \boldsymbol{U}_{1:k.1:k}^{-1}\right\Vert _{1}$\hfill{}\{update
inverse norms\}

\STATE $d_{k}\leftarrow a_{p_{k}q_{k}}-\boldsymbol{\ell}_{k,1:k-1}\boldsymbol{D}_{1:k-1,1:k-1}\boldsymbol{u}_{1:k-1,k}$\hfill{}\{recompute
diagonal entry due to pivoting\}

\STATE $n_{1}\leftarrow n_{1}-1$\hfill{}\{decrease leading block
for factorization\}

\ENDWHILE

\STATE $\boldsymbol{\ell}_{k+1:n,k}\leftarrow\frac{1}{d_{k}}\left(\hat{\vec{a}}_{\boldsymbol{p}_{k+1:n},q_{k}}-\vec{L}_{k+1:n,1:k-1}\boldsymbol{D}_{1:k-1,1:k-1}\vec{u}_{1:k-1,k}\right)$\label{line:ilu:crout-ell}\hfill{}\{fan-in
update $k$th column in $\boldsymbol{L}$\}

\STATE $\boldsymbol{u}_{k,k+1:n}\leftarrow\frac{1}{d_{k}}\left(\hat{\vec{a}}_{p_{k},\boldsymbol{q}_{k+1:n}}-\vec{\ell}_{k,1:k-1}\boldsymbol{D}_{1:k-1,1:k-1}\vec{U}_{1:k-1:k+1:n}\right)$\label{line:ilu:crout-u}\hfill{}\{fan-in
update $k$th row in $\boldsymbol{U}$\}

\STATE select largest $\left\lceil \alpha_{L}\text{nc}_{q_{k}}\right\rceil $
and $\left\lceil \alpha_{U}\text{nr}_{p_{k}}\right\rceil $ entries
in $\boldsymbol{\ell}_{k+1:n,k}$ and $\boldsymbol{u}_{k,k+1:n}$,
respectively \label{line:ilu:select-ell-u}\hfill{}\{scalability-oriented
dropping\}

\STATE drop entries $\left\{ \ell_{ik}\in\boldsymbol{\ell}_{k+1:n,k}\left|\kappa_{D}\tilde{\kappa}_{L}\vert\ell_{ik}\vert\leq\tau_{L}\right.\right\} $
and $\left\{ u_{kj}\in\boldsymbol{u}_{k,k+1:n}\left|\kappa_{D}\tilde{\kappa}_{U}\vert u_{kj}\vert\leq\tau_{U}\right.\right\} $
\label{line:ilu:num-drop}\hfill{}\{inverse-based dropping\}

\ENDFOR

\STATE $\vec{L}_{B}\leftarrow\vec{L}_{1:n_{1},1:n_{1}}$; $\boldsymbol{d}_{B}\leftarrow\vec{d}_{1:n_{1}}$;
$\vec{U}_{B}\leftarrow\vec{U}_{1:n_{1},1:n_{1}}$; $\vec{L}_{E}\leftarrow\vec{L}_{n_{1}+1:n,1:n_{1}}$;
$\vec{U}_{F}\leftarrow\vec{U}_{1:n_{1},n_{1}+1:n}$

\RETURN $\vec{L}_{B}$, $\boldsymbol{d}_{B}$, $\vec{U}_{B}$, $\vec{L}_{E}$,
$\vec{U}_{F}$, $\vec{p}$, $\vec{q}$, $n_{1}$

\end{algorithmic}
\end{algorithm}

\begin{figure}
\includegraphics[width=1\columnwidth]{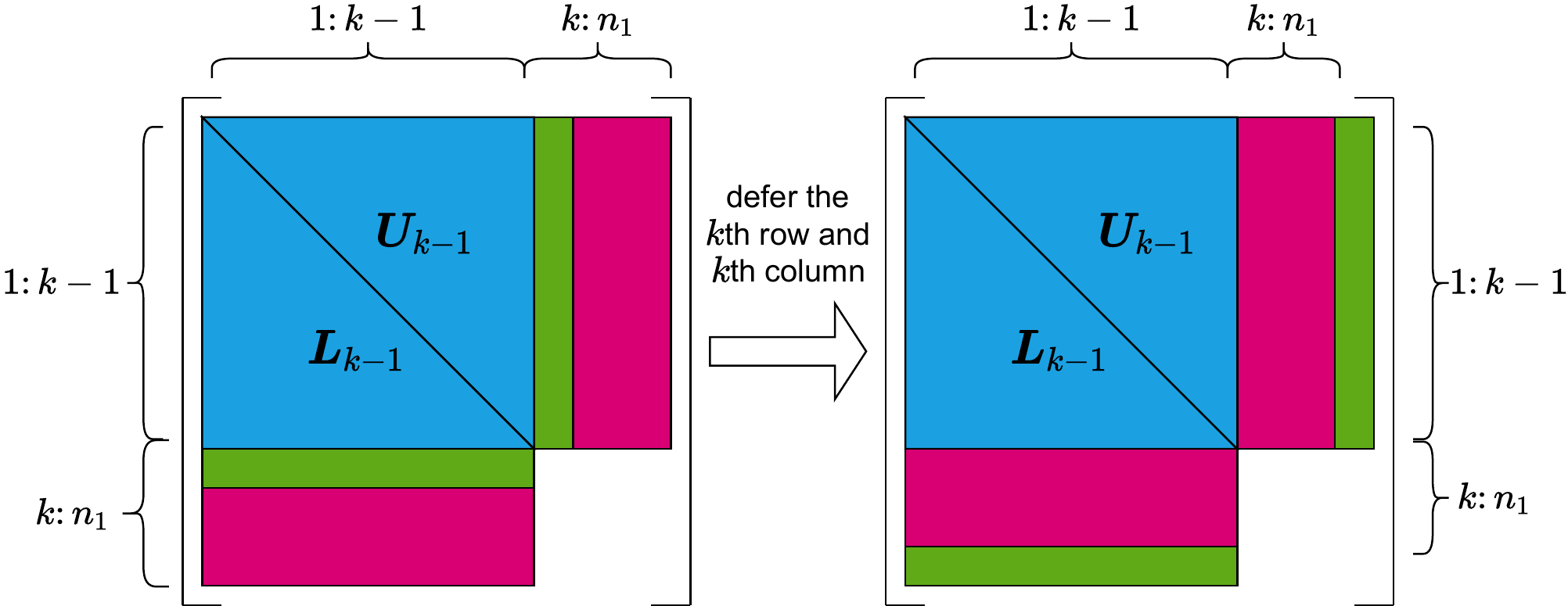}\caption{\label{fig:defer}Illustration of dynamic deferring. $\boldsymbol{L}_{k-1}$
and $\boldsymbol{U}_{k-1}$ are shorthand for $\boldsymbol{L}_{1:k-1,1:k-1}$
and $\boldsymbol{U}_{1:k-1,1:k-1}$, respectively. When encountering
ill-conditioned factors ($\left|d_{k}\right|<1/\kappa_{D}$ or $\max\left\{ \left\Vert \boldsymbol{L}_{1:k,1:k}^{-1}\right\Vert _{\infty},\left\Vert \boldsymbol{U}_{1:k,1:k}^{-1}\right\Vert _{1}\right\} >\kappa$)
due to $\boldsymbol{\ell}_{k,1:k-1}$ or $\boldsymbol{u}_{1:k-1,k}$
(green regions in left panel), we dynamically defer both $\boldsymbol{\ell}_{k,1:k-1}$
and $\boldsymbol{u}_{1:k-1,k}$ to next level (right panel).}
\end{figure}

Another core component in \textbf{ilu\_factorize} is the inverse-based
dynamic deferring. In particular, during the fan-in ILU, we dynamically
defer $\boldsymbol{\ell}_{k,1:k-1}$ and $\boldsymbol{u}_{1:k-1,k}$
(line~\ref{line:ilu:dyn-defer}) if we encounter small $\vert d_{k}\vert$
or large $\left\Vert \boldsymbol{L}_{1:k,1:k}^{-1}\right\Vert _{\infty}$
and $\left\Vert \boldsymbol{U}_{1:k,1:k}^{-1}\right\Vert _{1}$ (line~\ref{line:ilu:poor-factor}).
Figure~\ref{fig:defer} illustrates the process of dynamic deferring.
This deferring is similar to that in \cite{bollhofer2006multilevel,bollhofer2011ilupack};
in Section~\ref{subsec:Inverse-based-rook-pivoting}, we will extend
it to support rook pivoting. Note that the inverse-based dropping
in line~\ref{line:ilu:num-drop} can be replaced by a different dropping
strategy, such as that in \cite{mayer2006alternative}, but we utilized
the inverse-based dropping since we are already estimating the inverse
norms for deferring.

We note some implementation details. First, in line~\ref{line:ilu:select-ell-u},
we use quickselect \cite{hoare1961algorithm}, which has expected
linear time complexity. Second, since we need to access both rows
and columns of $\hat{\boldsymbol{A}}$ while computing the fan-in
updates (lines~\ref{line:ilu:crout-ell} and~\ref{line:ilu:crout-u}),
we need to store $\hat{\boldsymbol{A}}$ (or more precisely, $\boldsymbol{A}$)
in both row and column majors. We will describe the data structures
in Section~\ref{subsec:data-structures} and Appendix~\ref{sec:Flexible-array-based}.
Third, the output $\boldsymbol{L}_{E}$ and $\boldsymbol{U}_{F}$
will only be used to compute the Schur complement in Section~\ref{subsec:Schur-complement}.
Afterwards, $\boldsymbol{L}_{E}$ and $\boldsymbol{U}_{F}$ are discarded,
since they can be reconstructed from $\tilde{\boldsymbol{E}}=\hat{\boldsymbol{A}}_{\boldsymbol{p}_{n_{1}+1:n},\boldsymbol{q}_{1:n_{1}}}$
and $\tilde{\boldsymbol{F}}=\hat{\boldsymbol{A}}_{\boldsymbol{p}_{1:n_{1}},\boldsymbol{q}_{n_{1}+1:n}}$,
respectively as in (\ref{eq:two-lvl-ilu}).

\subsection{Inverse-based rook pivoting for coarse levels\label{subsec:Inverse-based-rook-pivoting}}

Dynamic deferring symmetrically permutes rows and columns. Such a
permutation strategy works well for reasonably well-conditioned matrices.
However, we observe that symmetric permutations alone sometimes lead
to relatively large Schur complements for highly ill-conditioned unsymmetric
systems. To overcome this issue, we introduce an \emph{inverse-based
rook pivoting} (\emph{IBRP}) for the fan-in ILU, by adapting the rook
pivoting \cite{poole2000rook} for complete LU factorization.

\begin{figure}
\centering{}\includegraphics[width=0.8\columnwidth]{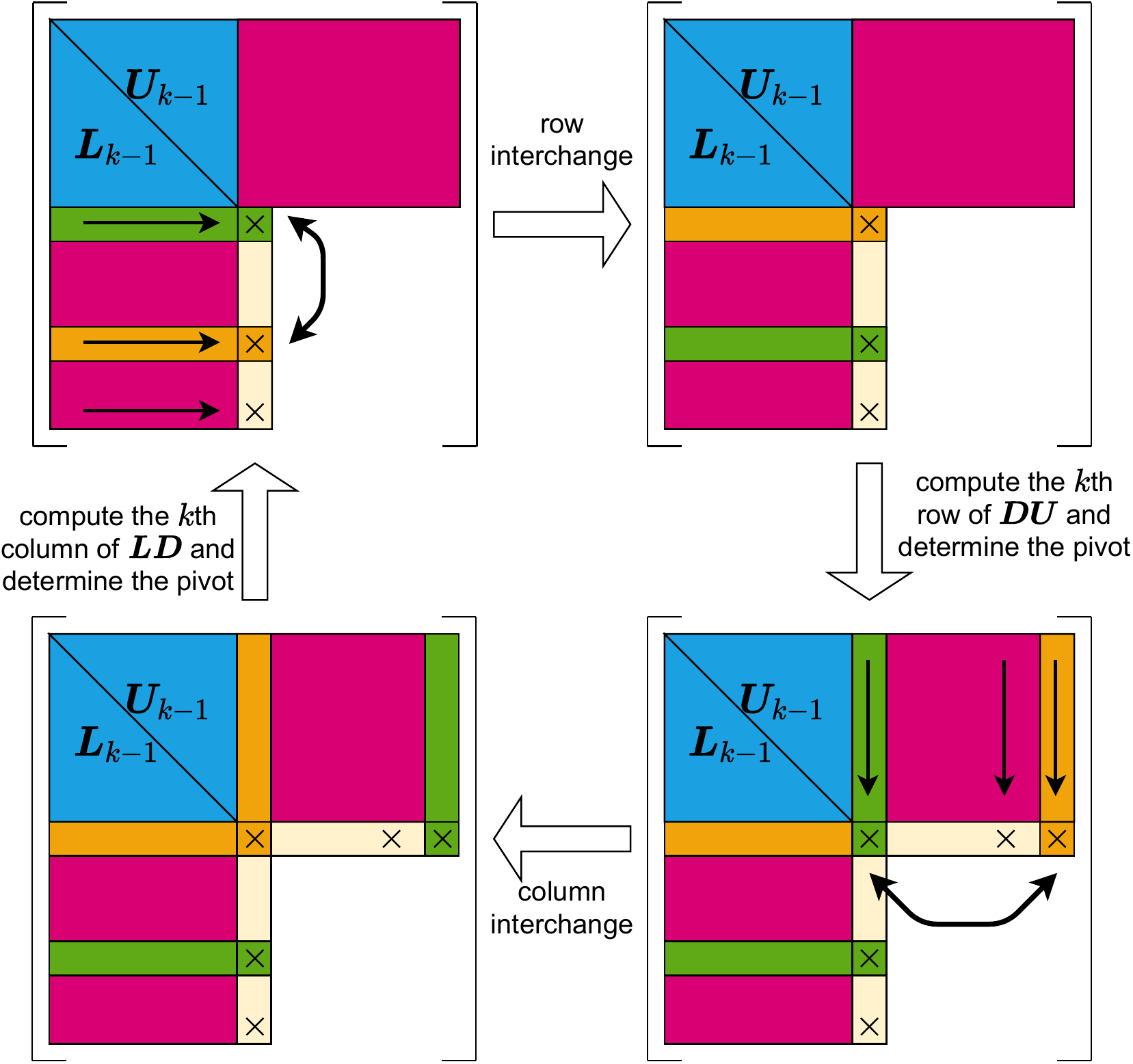}\caption{\label{fig:rook-pivot}Illustration of inverse-based rook pivoting
in HIF. At the upper-left panel, we update $d_{k}$ and nonzeros in
$\boldsymbol{\ell}_{k+1:n,k}$, find a pivot, and then interchange
$\boldsymbol{\ell}_{k,1:k-1}$ with the pivot row (indicated by curved
arrow). The lower-left panel shows the process for $\boldsymbol{U}$.}
\end{figure}

In the standard rook pivoting \cite{poole2000rook}, a pivot is found
by searching in the row and column in alternating order until its
magnitude is no smaller than those of all other entries in the row
and column within the Schur complement. This strategy has a comparable
cost as partial pivoting for dense matrices but enables superior stability.
However, in the context of fan-in ILU, only the $k$th row and column
of the Schur complement $\boldsymbol{S}$ are updated at the $k$th
step; the remaining part of the Schur complement is not available.
Hence, we must modify the pivoting procedure to interleave the search
with dynamic permutation and ``just-in-time'' fan-in updates. Figure~\ref{fig:rook-pivot}
illustrates one step of the IBRP. Note that this permutation is more
general than the dynamic deferring in Figure~\ref{fig:defer}, in
that it can exchange $\boldsymbol{\ell}_{k,1:k-1}$ with a row in
the middle of $\boldsymbol{L}_{k+1:n_{1},1:k-1}$, and similarly for
the rows in $\boldsymbol{U}$. Furthermore, the row and column interchanges
are not symmetric in general. Note that IBRP requires a more sophisticated
data structure to support the row and column interchanges, which we
will address in Section~\ref{subsec:data-structures}.

Besides the difference dictated by the fan-in update, there are two
other significant differences between IBRP and the standard rook pivoting.
First, we do not simply use the magnitude to determine the pivot,
since it may conflict with dynamic deferring. Instead, we add an inverse-based
constraint when searching the pivot, so that the pivot row in $\boldsymbol{L}$
and pivot column in $\boldsymbol{U}$ would not arbitrarily enlarge
the condition numbers of the $\boldsymbol{L}$ and $\boldsymbol{U}$,
respectively. Second, we do not locate the optimal rook pivot whose
magnitude is the largest among its row and column; instead, we impose
a maximum number of steps of IBRP, controlled by the parameter max\_steps.
For completeness, Algorithm~\ref{alg:rook_pivot} details the inverse-based
rook pivoting. Note that we do not scale $\hat{\boldsymbol{\ell}}$
and $\hat{\boldsymbol{u}}$ by $1/d_{k}$ in lines~\ref{line:pivot:ell}
and \ref{line:pivot:u} in Algorithm~\ref{alg:rook_pivot}, compared
to $\boldsymbol{\ell}_{k:n,k}$ and $\boldsymbol{u}_{k,k:n}$ in lines~\ref{line:ilu:crout-ell}
and \ref{line:ilu:crout-u} in Algorithm~\ref{alg:ilu_factor}. This
omission of scaling is for efficiency purposes, because the choice
of pivot does not depend on the scaling of the entries. 

\begin{algorithm}
\begin{raggedright}
\caption{\label{alg:rook_pivot}\textbf{ib\_rook\_pivot}$\left(\hat{\boldsymbol{A}},k,\boldsymbol{p},\boldsymbol{L},\boldsymbol{d},\boldsymbol{U},\boldsymbol{q},\kappa,n_{1},\text{max\_steps}\right)$}
\par\end{raggedright}
\begin{raggedright}
\textbf{inputs}:
\par\end{raggedright}
\begin{raggedright}
\hspace{1cm}\textbf{ }$\hat{\boldsymbol{A}}$: input scaled matrix
(i.e., $\hat{\boldsymbol{A}}=\boldsymbol{W}\boldsymbol{A}\boldsymbol{V}$,
passed in as $\boldsymbol{A}$, $\boldsymbol{W}$, and $\boldsymbol{V}$
separately)
\par\end{raggedright}
\begin{raggedright}
\hspace{1cm} $k$: step count in ILU factorization
\par\end{raggedright}
\begin{raggedright}
\hspace{1cm} $\boldsymbol{p},\boldsymbol{q}$: row and column permutation
vectors of $\boldsymbol{A}$, respectively
\par\end{raggedright}
\begin{raggedright}
\hspace{1cm} $\boldsymbol{L},\boldsymbol{U}$: $\boldsymbol{L}$
and $\boldsymbol{U}$ factors at step $k$, i.e, $\boldsymbol{L}_{B}\cup\boldsymbol{L}_{E}$
and $\boldsymbol{U}_{B}\cup\boldsymbol{U}{}_{F}$, respectively
\par\end{raggedright}
\begin{raggedright}
\hspace{1cm} $\boldsymbol{d}$: diagonal entries at step $k$
\par\end{raggedright}
\begin{raggedright}
\hspace{1cm} $\kappa$: inverse-norm threshold
\par\end{raggedright}
\begin{raggedright}
\hspace{1cm} $n_{1}$: leading block dimension
\par\end{raggedright}
\begin{raggedright}
\hspace{1cm} max\_steps: maximum number of rook pivoting steps
\par\end{raggedright}
\begin{raggedright}
\textbf{outputs}:
\par\end{raggedright}
\begin{raggedright}
\hspace{1cm} $\boldsymbol{p},\boldsymbol{q}$: updated row and column
permutation vectors, respectively
\par\end{raggedright}
\begin{raggedright}
\hspace{1cm} $\boldsymbol{L},\boldsymbol{U}$: updated $\boldsymbol{L}$
and $\boldsymbol{U}$ factors, respectively
\par\end{raggedright}
\begin{raggedright}
\hspace{1cm} $\boldsymbol{d}$: diagonal entries with updated $d_{k}$
\par\end{raggedright}
\begin{algorithmic}[1]

\FOR{$i=1$ to max\_steps}

\STATE $\hat{\boldsymbol{\ell}}_{k:n}\leftarrow\hat{\vec{a}}_{\boldsymbol{p}_{k:n},q_{k}}-\vec{L}_{k:n,1:k-1}\vec{D}_{1:k-1,1:k-1}\vec{u}_{1:k-1,k}$\label{line:pivot:ell}\hfill{}\{fan-in
update of $d_{k}$ and $\boldsymbol{\ell}_{k+1:n,k}$ without scaling
by diagonal\}

\STATE $r\leftarrow\arg\max$$\left\{ \text{\ensuremath{\left|\hat{\ell}_{r}\right|}}\left|k+1\leq r\le n_{1}\text{ and }\left\Vert \begin{bmatrix}\boldsymbol{L}_{1:k-1,1:k-1} & \boldsymbol{0}\\
\boldsymbol{L}_{r,1:k-1} & \ell_{r,k}
\end{bmatrix}^{-1}\right\Vert _{\infty}\le\kappa\right.\right\} $ \label{line:pivot:pvtell}\hfill{}\{find pivot in $\hat{\boldsymbol{\ell}}_{k+1:n}$\}

\IF{$r\neq\emptyset$ \textbf{and} $\left|\hat{\ell}_{k}\right|<\left|\hat{\ell}_{r}\right|$}

\STATE $\boldsymbol{\ell}_{k,1:k-1}\leftrightarrow\boldsymbol{\ell}_{r,1:k-1}$;
$p_{k}\leftrightarrow p_{r}$; $\hat{\ell}_{k}\leftrightarrow\hat{\ell}_{r}$\label{line:pivot:swapell}\hfill{}\{perform
row interchange\}

\ELSIF{$i>1$} 

\STATE $d_{k}\leftarrow\hat{\ell}_{k}$\textbf{; break}\hfill{}\{extract
$d_{k}$ from $\hat{\boldsymbol{\ell}}_{k:n}$\}

\ENDIF

\STATE $\hat{\boldsymbol{u}}_{k:n}\leftarrow\hat{\vec{a}}_{p_{k},\boldsymbol{q}_{k:n}}-\vec{\ell}_{k,1:k-1}\vec{D}_{1:k-1,1:k-1}\vec{U}_{1:k-1:k:n}$\label{line:pivot:u}\hfill{}\{fan-in
update of $d_{k}$ and $\boldsymbol{u}_{k,k+1:n}$ without scaling
by diagonal\}

\STATE $c\leftarrow\arg\max$$\left\{ \text{\ensuremath{\left|\hat{u}_{c}\right|}}\left|k+1\leq c\le n_{1}\text{ and }\left\Vert \begin{bmatrix}\boldsymbol{U}_{1:k-1,1:k-1} & \boldsymbol{U}_{1:k-1,c}\\
\boldsymbol{0} & u_{k,c}
\end{bmatrix}^{-1}\right\Vert _{1}\le\kappa\right.\right\} $ \label{line:pivot:pvtu}\hfill{}\{find pivot in $\hat{\boldsymbol{u}}_{k+1:n}$\}

\IF{ $c\neq\emptyset$ \textbf{and} $\left|\hat{u}_{k}\right|<\left|\hat{u}_{c}\right|$}

\STATE $\boldsymbol{u}_{1:k-1,k}\leftrightarrow\boldsymbol{u}_{1:k-1,c}$;
$q_{k}\leftrightarrow q_{c}$; $\hat{u}_{k}\leftrightarrow\hat{u}_{r}$\label{line:pivot:swapu}\hfill{}\{perform
column interchange\}

\ELSE

\STATE $d_{k}\leftarrow\hat{u}_{k}$\textbf{; break}\hfill{}\{extract
$d_{k}$ from $\hat{\boldsymbol{u}}_{k:n}$\}

\ENDIF

\ENDFOR

\RETURN $\boldsymbol{p},\ \boldsymbol{L},\ \boldsymbol{d},\ \boldsymbol{U},\ \boldsymbol{q}$

\end{algorithmic}
\end{algorithm}

We note an important practical issue. The IBRP can result in significantly
denser $\boldsymbol{L}$ and $\boldsymbol{U}$ factors because the
pivoting may undo the effects of the fill-reduction reordering in
the preprocessing step. Hence, we enable IBRP only on the coarser
levels (typically for $\text{level}>1$) when dynamic deferring is
found to be ineffective. At the coarser levels, we also enlarge the
scalability-oriented fill factors $\alpha_{L}$ and $\alpha_{U}$
to preserve more \emph{fills} introduced by IBRP. Our experiments
show that applying IBRP on coarser levels can significantly reduce
the size of the final Schur complement for some challenging problems,
as demonstrated by the example in Table~\ref{tab:effect-rook-pivot}.

\begin{table}
\caption{\label{tab:effect-rook-pivot}Effectiveness of IBRP for reducing the
final Schur complement size for the testing matrix \textsf{shyy161}
from the SuiteSparse Matrix Collection \cite{davis2011university}.
$\text{dim}(\mathcal{N})$ indicates the dimension of the null space,
and $\boldsymbol{S}_{m}$ denotes the final Schur complement in HIF.
We factorized the testing matrix using HIF with and without IBRP,
and we report the factorization time (factor time) and nnz ratio (i.e.,
$\nicefrac{\text{nnz}(\boldsymbol{M})}{\text{nnz}(\boldsymbol{A})}$,
aka fill ratio in \cite{li2011supernodal}). We solved the consistent
system with right-hand side $\boldsymbol{b}=\boldsymbol{A}\boldsymbol{1}$
(where $\boldsymbol{1}=[1,1,\dots,1]^{T}$) using GMRES(30) with relative
residual tolerance $10^{-6}$. The computing environment for running
the test can be found in Section~\ref{sec:Illustration-of-HIFIR}.
Timing results for GMRES were negligible thus omitted. }

\centering{}%
\begin{tabular}{cccc||cccccc}
\hline 
\multirow{2}{*}{$n$} & \multirow{2}{*}{nnz} & \multirow{2}{*}{$\text{dim}(\mathcal{N}(\boldsymbol{A}))$} & \multirow{2}{*}{w/ IBRP} & \multicolumn{2}{c}{$\boldsymbol{S}_{m}$ in HIF} & fac. & \multirow{2}{*}{\#levels} & nnz & GMRES\tabularnewline
 &  &  &  & $n$ & $\text{dim}(\mathcal{N}(\boldsymbol{S}_{m}))$ & time &  & ratio & iter.\tabularnewline
\hline 
\hline 
\multirow{2}{*}{76,480} & \multirow{2}{*}{329,762} & \multirow{2}{*}{$\gtrsim50$} & yes & 1,407 & 50 & 2.51 & 5 & 12.0 & 15\tabularnewline
 &  &  & no & 9,394 & 48 & 226 & 2 & 304 & 2\tabularnewline
\hline 
\end{tabular}
\end{table}

\subsection{Computing and factorizing Schur complements\label{subsec:Schur-complement}}

After finishing \textbf{ilu\_factorize}, we need to compute the Schur
complement $\boldsymbol{S}$ based on (\ref{eq:schur-complement}).
This step involves a sparse matrix-matrix (SpMM) multiplication, for
which we adopt the algorithm as described in \cite{bank1993sparse}.
The space and time complexity of SpMM depends on the nonzeros in $\boldsymbol{L}_{E}$
and $\boldsymbol{U}_{F}$. To achieve near-linear complexity, we apply
scalability-oriented dropping \textit{before} SpMM to the rows and
columns of $\boldsymbol{L}_{E}$ and $\boldsymbol{U}_{F}$, respectively.
Recall that in \textbf{ilu\_factorize} we already applied scalability-oriented
dropping to the columns and rows of $\boldsymbol{L}_{E}$ and $\boldsymbol{U}_{F}$,
respectively. Hence, the nonzeros both rows and columns in $\boldsymbol{L}_{E}$
and in $\boldsymbol{U}_{F}$ are well controlled, allowing us to effectively
bound the complexity of $\boldsymbol{S}$. In contrast, if we applied
dropping \textit{after} SpMM, the complexity of SpMM may be higher.
Note that after computing $\boldsymbol{S}$, we drop both $\boldsymbol{L}_{E}$
and $\boldsymbol{U}_{F}$ as they can be reconstructed from the other
terms as in (\ref{eq:two-lvl-ilu}), and then factorize $\boldsymbol{S}$
recursively.

In HIF, the final Schur complement is factorized by rank-revealing
QR (or truncated QRCP) in order to guarantee the stability of the
preconditioner. Algorithm~\ref{alg:dense_factor} outlines the procedure.
We note a couple of details in the algorithm. First, line~\ref{line:dense_factor:qrcp}
computes QRCP $\boldsymbol{S}\boldsymbol{P}=\boldsymbol{Q}\boldsymbol{R}$
for $\boldsymbol{S}\in\mathbb{R}^{n_{S}\times n_{S}}$, where $\boldsymbol{P}$
is a permutation matrix, and $r_{11}\ge r_{22}\ge\ldots\ge r_{n_{S}n_{S}}\geq0$.
The QRCP has a time complexity of $\mathcal{O}\left(n_{S}^{3}\right)$.
Second, lines~\ref{line:dense_factor:kappa-start}--\ref{line:dense_factor:kappa-end}
determine the numerical rank of $\boldsymbol{S}$ by comparing the
estimated condition number against a threshold $\kappa_{\text{rrqr}}$,
which defaults to $\kappa_{\text{rrqr}}=\epsilon_{\text{mach}}^{-2/3}$.
We estimate the $2$-norm condition number using the incremental estimator
in \cite{bischof1990incremental}, which has a linear complexity per
step. Hence, the overall computational cost is dominated by QRCP.
For efficiency, we implement the QRCP and condition-number estimator
using LAPACK kernel functions \textsf{xGEQP3 }and \textsf{xLAIC1},
respectively.

\begin{algorithm}
\caption{\label{alg:dense_factor}\textbf{rrqr\_final\_schur}$\left(\boldsymbol{S},\kappa_{\text{rrqr}}\right)$}

\begin{raggedright}
\textbf{inputs}:
\par\end{raggedright}
\begin{raggedright}
\hspace{1cm} $\boldsymbol{S}$: final Schur complement of size $n_{s}\times n_{s}$
\par\end{raggedright}
\begin{raggedright}
\hspace{1cm} $\kappa_{\text{rrqr}}$: condition number threshold
for determining numerical rank (default value is $\epsilon_{\text{mach}}^{-2/3}$)
\par\end{raggedright}
\begin{raggedright}
\textbf{outputs}:
\par\end{raggedright}
\begin{raggedright}
\hspace{1cm} \textbf{$\boldsymbol{Q},\boldsymbol{R}$}: $\boldsymbol{Q}$
and $\boldsymbol{R}$ factors of $\boldsymbol{S}$
\par\end{raggedright}
\begin{raggedright}
\hspace{1cm} $\boldsymbol{p}$: permutation vector in QRCP
\par\end{raggedright}
\begin{raggedright}
\hspace{1cm} $r$: numerical rank of $\boldsymbol{S}$
\par\end{raggedright}
\begin{algorithmic}[1]

\STATE \label{line:dense_factor:qrcp}factorize $\boldsymbol{S}=\boldsymbol{Q}\boldsymbol{R}\boldsymbol{P}^{T}$\hfill{}\{compute
QRCP via \textsf{xGEQP3}\}

\FOR{$r=1$ \TO $n_{S}$}\label{line:dense_factor:kappa-start} 

\STATE \label{line:dense_factor:est-kappa}update estimated $\left\Vert \boldsymbol{R}_{1:r,1:r}\right\Vert _{2}$
and $\left\Vert \boldsymbol{R}_{1:r,1:r}^{-1}\right\Vert _{2}$\hfill{}\{incremental
condition number estimation via \textsf{xLAIC1}\}

\IF{$\left\Vert \boldsymbol{R}_{1:r,1:r}\right\Vert _{2}\left\Vert \boldsymbol{R}_{1:r,1:r}^{-1}\right\Vert _{2}<\kappa_{\text{rrqr}}$}

\STATE $r\leftarrow r+1$

\STATE \textbf{break}

\ENDIF

\ENDFOR\label{line:dense_factor:kappa-end}

\RETURN $\boldsymbol{Q}$, $\boldsymbol{R}$, $\boldsymbol{p}$,
$r$

\end{algorithmic}
\end{algorithm}

Due to the cubic time complexity of QRCP, we would like to make $n_{S}$
as small as possible, and ideally have $n_{S}=\mathcal{O}\left(\sqrt[3]{n}\right)$,
where $n$ is the number of unknowns in the original system. This
high complexity was the motivation to use IBRP within multilevel ILU.
In addition, we would also like to prevent having too many low-quality
ILU levels. As a tradeoff, we trigger QRCP based on the following
criteria. In \textbf{ilu\_factorize}, let $n_{0}$ be the initial
leading block dimension after preprocessing, i.e., the input $n_{1}$
in Algorithm~\ref{alg:ilu_factor}, and let $d$ be the total number
of dynamic deferrals. If more than $75\%$ entries are dynamically
deferred during fan-in update, i.e., $\nicefrac{d}{n_{0}}\ge0.75$,
then we discard incomplete factorization in that level and use RRQR
on its input. In addition, if more than $60\%$ entries are dynamically
deferred, i.e., $\nicefrac{d}{n_{0}}\ge0.6$, then we apply RRQR on
the remainder Schur complement.

\subsection{Preprocessing with static deferring\label{subsec:Symmetric-and-asymmetric}}

In Figure~\ref{fig:hif-workflow}, an important step was the preprocessing
at each step, which computes a block structure $\begin{bmatrix}\hat{\boldsymbol{B}} & \hat{\boldsymbol{F}}\\
\hat{\boldsymbol{E}} & \hat{\boldsymbol{C}}
\end{bmatrix}$. The preprocessing in HIF performs equilibration, static deferring,
and fill-reduction reordering, in that order. Equilibration improves
stability by computing $\boldsymbol{P}^{T}\boldsymbol{W}\boldsymbol{A}\boldsymbol{V}\boldsymbol{Q}=\tilde{\boldsymbol{A}}$,
where $\boldsymbol{W}$ and $\boldsymbol{V}$ correspond to row and
column scaling matrices, and $\boldsymbol{P}$ and $\boldsymbol{Q}$
correspond to row and column permutation matrices. We utilize MC64
\cite{duff2001algorithms}, which computes unsymmetric equilibration,
for which $\boldsymbol{Q}=\boldsymbol{I}$. For (nearly) pattern symmetric
levels, we symmetrize the output of MC64 by setting $\boldsymbol{W}=\boldsymbol{V}=\sqrt{\boldsymbol{W}\boldsymbol{V}}$
and $\boldsymbol{Q}=\boldsymbol{P}$ as in \cite{hsl_mc64}. Note
that for structurally singular systems, MC64 sometimes yields unstable
scaling factors that may be arbitrarily large or small \cite{jiao2021approximate}.
The symmetrization process overcomes the issue; for unsymmetric equilibration,
we also symmetrize the scaling factors by setting $w_{i}=v_{i}=\sqrt{w_{i}v_{i}}$
if $\max\left\{ w_{i},v_{i}\right\} /\min\left\{ \text{\ensuremath{w_{i}},\ensuremath{v_{i}}}\right\} >\beta$,
where $\beta$ is 1000 by default. After symmetric equilibration,
we permute a row and its corresponding column to the lower-right corner
if its diagonal entry is nearly zero. We refer to this process as
\emph{static deferring}, which naturally yields a block structure
$\begin{bmatrix}\boldsymbol{B} & \boldsymbol{F}\\
\boldsymbol{E} & \hat{\boldsymbol{C}}
\end{bmatrix}$. Afterwards, we apply fill-reduction reordering on $\boldsymbol{B}$
and then permute $\boldsymbol{E}$ and $\boldsymbol{F}$ correspondingly,
which leads to the final block structure $\begin{bmatrix}\hat{\boldsymbol{B}} & \hat{\boldsymbol{F}}\\
\hat{\boldsymbol{E}} & \hat{\boldsymbol{C}}
\end{bmatrix}$. We use RCM \cite{george1971computer} and AMD \cite{amestoy2004algorithm}
for symmetric and unsymmetric reordering, respectively, because RCM
is more efficient and has been shown to work better for symmetric
ILU \cite{benzi1999orderings,gupta2010adaptive}. 

\subsection{\label{subsec:data-structures}Efficient data structures and complexity
analysis}

To implement HIF efficiently, we must perform all its core operations
in linear time with respect to the number of nonzeros. In particular,
the fan-in updates of $\boldsymbol{\ell}_{k+1:n,k}$ and $\boldsymbol{u}_{k+1:n,k}$
at $k$th \textbf{ilu\_factorize} (cf. lines~\ref{line:ilu:crout-ell}
and \ref{line:ilu:crout-u} in Algorithm~\ref{alg:ilu_factor} and
lines~\ref{line:pivot:ell} and \ref{line:pivot:u} in Algorithm~\ref{alg:rook_pivot})
requires efficient sequential access of $\boldsymbol{L}$ and $\boldsymbol{U}$
in both rows and columns. More importantly, the deferring and pivoting
operations require interchanging rows and columns in $\boldsymbol{L}$
and $\boldsymbol{U}$ dynamically. Although the data structure in
\cite{li2003crout} supports fan-in updates efficiently, it does not
support deferring or rook pivoting. Under these considerations, we
developed flexible, three-tiered data structures, which extended the
data structure in \cite{li2003crout} to support deferring and pivoting,
as we describe in Appendix~\ref{sec:Flexible-array-based}. This
three-tiered data structure augments the standard compressed sparse
column or row (aka CSC and CSR) formats either partially or fully.
We use the partially augmented version when rook pivoting is disabled,
especially at the top levels, since it has a smaller memory footprint;
the fully augmented version is used when rook pivoting is enabled,
which in general occurs only at coarser levels. 

Assuming the number of nonzeros per row and column in the input mesh
is bounded by a constant, which is typically the case for linear systems
arising from partial differential equations, the augmented data structures
enable all the core components in HIF to be performed in linear time
at each level. We outline the analysis as follows. First, the fan-in
update can be performed in linear time proportional to the number
of nonzeros, thanks to the use of augmented data structure \cite{li2003crout}.
Second, the incremental update of $\left\Vert \boldsymbol{L}_{1:k,1:k}^{-1}\right\Vert _{\infty}$
and $\left\Vert \boldsymbol{U}_{1:k,1:k}^{-1}\right\Vert _{1}$ also
costs linear time with respect to the number of nonzeros in $\boldsymbol{L}$
and $\boldsymbol{U}$, respectively (see e.g., \cite{Golub13MC}).
Third, the partially augmented data structures allow the permutations
in deferring can be performed in linear time complexity, as shown
in \cite{chen2021hilucsi}. Fourth, the fully augmented data structures
allow each row interchange in rook pivoting to be performed in time
proportional to the number of nonzeros, assuming the number of nonzeros
in each row is a constant. This step is the most complicated, and
we defer its analysis to Appendix~\ref{sec:Time-complexity-rook-pivoting}.
Fifth, the time complexity of SpMM in Section~\ref{subsec:Schur-complement}
is proportional to the number of nonzeros in the product. Since the
scalability-oriented dropping ensures that each row and column of
$\boldsymbol{L}_{E}$ and $\boldsymbol{U}_{F}$ is bounded by a constant
$C$, so the number of nonzeros per row and per column in the product
$\boldsymbol{L}_{E}\boldsymbol{D}_{B}\boldsymbol{U}_{F}$ is bounded
by $C^{2}$. Finally, assuming the number of rows in the final Schur
complement is bounded by $\mathcal{O}(n^{1/3})$, we conclude that
HIF guarantees linear time complexity per level (excluding its preprocessing
steps) under the assumptions as mentioned above.

Note that HIF does not guarantee that the number of levels is bounded
by a constant. Furthermore, the time complexity of AMD reordering
is quadratic in the worst-case case \cite{Heggernes01CCM}. Hence,
the total time complexity of HIF may not be linear. Nevertheless,
the number of levels is typically a small constant, and the computational
cost of AMD is typically negligible. Hence, we do observe linear asymptotic
growth for the overall HIF across all levels empirically.

\section{\label{sec:Applying-HIFIR}Multilevel triangular solves and matrix-vector
multiplications}

To use the HIF as a preconditioner, we typically need a procedure
similar to triangular solves. We shall refer to it as a \emph{multilevel
triangular solver}. As an illustration, let $\boldsymbol{y}=\begin{bmatrix}\boldsymbol{y}_{1}\\
\boldsymbol{y}_{2}
\end{bmatrix}$ be a block vector corresponding to $\tilde{\boldsymbol{B}}$ and
$\boldsymbol{C}$ in the two-level ILU $\boldsymbol{M}$ in (\ref{eq:two-lvl-ilu}).
In this case, $\boldsymbol{W}^{-1}\boldsymbol{P}\boldsymbol{M}\boldsymbol{Q}^{T}\boldsymbol{V}^{-1}$
is a preconditioner of $\boldsymbol{A}$. The multilevel solver computes
\begin{equation}
\boldsymbol{M}^{g}\boldsymbol{y}=\begin{bmatrix}\tilde{\boldsymbol{B}}^{-1}\boldsymbol{y}_{1}\\
\boldsymbol{0}
\end{bmatrix}+\begin{bmatrix}-\tilde{\boldsymbol{B}}^{-1}\boldsymbol{F}\\
\boldsymbol{I}
\end{bmatrix}\boldsymbol{S}^{g}\left(\boldsymbol{y}_{2}-\boldsymbol{E}\tilde{\boldsymbol{B}}^{-1}\boldsymbol{y}_{1}\right).\label{eq:back-solve}
\end{equation}
We further compute $\boldsymbol{S}^{g}\boldsymbol{y}_{1}$ recursively,
which leads to the multilevel solver. Algorithm~\ref{alg:hif_solve}
details this recursive procedure. Note that Algorithm~\ref{alg:hif_solve}
also supports the use of $\boldsymbol{M}^{H}$ to construct a preconditioner
for $\boldsymbol{A}^{H}$ based on (\ref{eq:conjugate-trans}).

\begin{algorithm}
\caption{\label{alg:hif_solve}\textbf{hif\_solve}$\left(\text{\ensuremath{\boldsymbol{M}}},\boldsymbol{y},r,\text{trans},\text{level}=1\right)$}

\begin{raggedright}
\textbf{inputs}:
\par\end{raggedright}
\begin{raggedright}
\hspace{1cm} $\text{\ensuremath{\boldsymbol{M}}}$: a structure containing
HIF preconditioner $\boldsymbol{M}$
\par\end{raggedright}
\begin{raggedright}
\hspace{1cm} $\boldsymbol{y}$: right-hand side vector
\par\end{raggedright}
\begin{raggedright}
\hspace{1cm} $r$: numerical rank used in the final RRQR factorization
\par\end{raggedright}
\begin{raggedright}
\hspace{1cm} trans: Boolean flag indicating whether to solve $\boldsymbol{M}^{gH}\boldsymbol{y}$
\par\end{raggedright}
\begin{raggedright}
\hspace{1cm} level: level counter in the HIF preconditioner (default
to 1 at top level)
\par\end{raggedright}
\begin{raggedright}
\textbf{output}:
\par\end{raggedright}
\begin{raggedright}
\hspace{1cm} $\boldsymbol{x}$: solution vector $\boldsymbol{M}^{g}\boldsymbol{y}$
or $\boldsymbol{M}^{gH}\boldsymbol{y}$
\par\end{raggedright}
\begin{algorithmic}[1]

\STATE extract $n_{1}$, $n$, $\boldsymbol{L}_{B}$, $\boldsymbol{D}_{B}$,
$\boldsymbol{U}_{B}$, $\boldsymbol{E}$, $\boldsymbol{F}$, $\boldsymbol{P}$,
$\boldsymbol{Q}$, $\boldsymbol{W}$, and $\boldsymbol{V}$ for current
level from $\boldsymbol{M}$

\IF{trans}

\STATE $\boldsymbol{L}_{B},\boldsymbol{U}_{B}\leftarrow\boldsymbol{U}_{B}^{H},\boldsymbol{L}_{B}^{H}$;
$\boldsymbol{D}_{B}\leftarrow\bar{\boldsymbol{D}}_{B}$\hfill{}\{construct
conjugate transpose\}

\STATE $\boldsymbol{E},\boldsymbol{F}\leftarrow\boldsymbol{F}^{H},\boldsymbol{E}^{H}$

\STATE $\boldsymbol{W}\leftrightarrow\boldsymbol{V}$; $\boldsymbol{P}\leftrightarrow\boldsymbol{Q}$

\ENDIF

\STATE $\boldsymbol{y}\leftarrow\boldsymbol{P}^{T}\boldsymbol{W}\boldsymbol{y}$\hfill{}\{scale
and permute right-hand side\}

\STATE $\boldsymbol{x}_{1:n_{1}}\leftarrow\boldsymbol{U}_{B}^{-1}\boldsymbol{D}_{B}^{-1}\boldsymbol{L}_{B}^{-1}\boldsymbol{y}_{1:n_{1}}$

\STATE $\boldsymbol{x}_{n_{1}+1:n}\leftarrow\boldsymbol{y}_{n_{1}+1:n}-\boldsymbol{E}\boldsymbol{x}_{1:n_{1}}$

\IF{ final level}

\STATE extract $\hat{\boldsymbol{Q}}$, $\hat{\boldsymbol{R}}$,
and $\hat{\boldsymbol{P}}$ of QRCP of $\boldsymbol{S}_{m}$ (or of
$\boldsymbol{S}_{m}^{H}$ if trans) from $\boldsymbol{M}$

\STATE $\boldsymbol{x}_{n_{1}+1:n}\leftarrow\hat{\boldsymbol{P}}_{:,1:r}\hat{\boldsymbol{R}}_{r}^{-1}\hat{\boldsymbol{Q}}_{:,1:r}^{H}\thinspace\boldsymbol{x}_{n_{1}+1:n}$\hfill{}\{solve
on final Schur complement with numerical rank $r$\}

\ELSE

\STATE $\boldsymbol{x}_{n_{1}+1:n}\leftarrow\text{\textbf{hif\_solve}}\left(\boldsymbol{M},\boldsymbol{x}_{n_{1}+1:n},r,\text{tran},\text{level}+1\right)$

\ENDIF

\STATE $\boldsymbol{x}_{1:n_{1}}\leftarrow\boldsymbol{y}_{1:n_{1}}-\boldsymbol{F}\boldsymbol{x}_{n_{1}+1:n}$

\STATE $\boldsymbol{x}_{1:n_{1}}\leftarrow\boldsymbol{U}_{B}^{-1}\boldsymbol{D}_{B}^{-1}\boldsymbol{L}_{B}^{-1}\boldsymbol{x}_{1:n_{1}}$

\RETURN $\boldsymbol{V}\boldsymbol{Q}\boldsymbol{x}$\hfill{}\{scale
and permute $\boldsymbol{x}$\}

\end{algorithmic}
\end{algorithm}

When solving inconsistent systems, such as the computation of the
null-space vector, using $\boldsymbol{G}=\tilde{\boldsymbol{M}}^{g}$
as the preconditioning operator may be insufficient. In this setting,
it is desirable to enable iterative refinement with HIF to construct
a variable preconditioner for FGMRES, as described in \cite{jiao2021approximate}.
Specifically, given equation $\boldsymbol{A}\boldsymbol{v}=\boldsymbol{q}$,
HIFIR iteratively computes 
\begin{equation}
\boldsymbol{v}_{j}=(\boldsymbol{I}-\boldsymbol{M}^{g}\boldsymbol{A})\boldsymbol{v}_{j-1}+\boldsymbol{M}^{g}\boldsymbol{q}=\boldsymbol{v}_{j-1}-\boldsymbol{M}^{g}\boldsymbol{r}_{j-1}\label{eq:iterative-refinement}
\end{equation}
for $j=1,2,\ldots$, where $\boldsymbol{r}_{j-1}=\boldsymbol{q}-\boldsymbol{A}\boldsymbol{v}_{j-1}$
is the residual vector with $\boldsymbol{v}_{0}=\boldsymbol{0}$,
and $\boldsymbol{M}^{g}\boldsymbol{r}_{j-1}$ is computed by Algorithm~\ref{alg:hif_solve}.
Note that when computing null-space vectors, the condition number
of $\boldsymbol{R}_{1}$ in the RRQR of $\boldsymbol{S}_{m}$ can
be as large as $1/\epsilon_{\text{mach}}$, so it would have a different
numerical rank of $\boldsymbol{S}_{m}$ compared to that when solving
consistent systems as described in Section~\ref{subsec:Schur-complement}.
For this reason, Algorithm~\ref{alg:hif_solve} has a parameter $r$
to allow the user passing in different numerical ranks of $\boldsymbol{S}_{m}$.
When choosing the number of iterations in the iterative refinement,
HIFIR terminates the iteration after a maximum number of iterations.
We increase this upper bound on the iterations at every restart for
FGMRES. We refer readers to \cite{jiao2021approximate} for the analysis
of this adaptive procedure.

Besides multilevel triangular solves, HIF also supports the computation
of the multiplication $\tilde{\boldsymbol{M}}$ with a vector. As
an illustration, given $\boldsymbol{x}=\begin{bmatrix}\boldsymbol{x}_{1}\\
\boldsymbol{x}_{2}
\end{bmatrix}$ in a two-level factorization, the multiplication can be written as
\begin{align}
\tilde{\boldsymbol{M}}\boldsymbol{x} & =\begin{bmatrix}\tilde{\boldsymbol{B}}\\
\boldsymbol{E}
\end{bmatrix}\boldsymbol{x}_{1}+\begin{bmatrix}\boldsymbol{I}\\
\boldsymbol{E}\tilde{\boldsymbol{B}}^{-1}
\end{bmatrix}\boldsymbol{F}\boldsymbol{x}_{2}+\begin{bmatrix}\boldsymbol{0}\\
\boldsymbol{S}\boldsymbol{x}_{2}
\end{bmatrix}.\label{eq:mat-vector}
\end{align}
The multiplication of $\boldsymbol{x}_{2}$ by the Schur complement,
i.e. $\boldsymbol{S}\boldsymbol{x}_{2}$, is then computed recursively,
leading to the multilevel matrix-vector multiplication. The control
flow of this recursive algorithm is similar to that of Algorithm~\ref{alg:hif_solve},
and it is helpful when we need $\tilde{\boldsymbol{M}}$ as an approximation
of $\boldsymbol{A}$ (instead of using $\tilde{\boldsymbol{M}}^{g}$
as an approximation to $\boldsymbol{A}^{g}$). In addition, it is
also useful for some advanced preconditioners for singular systems,
which we will report in the future.

\section{\label{sec:software-library}Software design and user interfaces}

In this section, we describe our software design and the user interfaces
of HIFIR, including its C++ programming interface as well as the high-level
interfaces for Python and MATLAB. We refer the readers to the official
documentation \url{https://hifirworks.github.io/hifir/} for more
detailed documentation of HIFIR.

\subsection{Design considerations}

When designing HIFIR, we have focused on three key factors: efficiency,
flexibility, and ease-of-integration into other codes. Under these
considerations, we chose to implement the core components of HIFIR
in C++-11 using C++ templates in a header-only fashion while providing
high-level interfaces in Python and MATLAB. C++ is a powerful programming
language for scientific computing and is highly efficient. For efficient
kernel computations, we link HIFIR with LAPACK and we ease the linking
by leveraging compiler directives and macros. In terms of flexibility,
HIFIR supports different data types, such as \textsf{float}, \textsf{double},
and \textsf{std::complex<double>}, through C++ templates. In addition,
mixed data types are supported, for example, to compute the factorization
in double precision while exporting it in single precision during
the solve step. In addition, the use of C++ also makes it easier to
integrate with other linear-algebra packages (such as Eigen \cite{eigenweb},
Blaze \cite{iglberger2012expression}, SuperLU \cite{li2005overview},
PETSc \cite{petsc-user-ref}, etc.) and scientific-computing packages
(such as Trilinos \cite{heroux2005overview}).  The use of C++ in
HIFIR eases its integration with such libraries. The header-only design
of HIFIR is similar to that of Eigen \cite{eigenweb}. This design
significantly simplifies the installation and build process in that
the user only needs to include the header file \textsf{hifir.hpp}
and link with LAPACK libraries. For ease of prototyping, we also provide
high-level interfaces in Python and MATLAB through\textsf{ hifir4py}
and \textsf{hifir4m}, respectively; see Section~\ref{subsec:interface-for-Python-MATLAB}.

\subsection{Control parameters\label{subsec:Control-parameters}}

The algorithm described in Section~\ref{sec:Algorithmic-components}
involved several control parameters, as we summarize some key parameters
in Figure~\ref{fig:parameters-struct} as a C++ structure. The first
six parameters correspond to $\alpha_{L}$, $\alpha_{U}$, $\kappa$,
$\kappa_{D}$, $\tau_{L}$, and $\tau_{U}$ for the top level in Algorithm~\ref{alg:ilu_factor};
note that HIF automatically adapt these parameters for coarser levels.
Their default values were obtained based on the experiments in \cite{chen2021hilucsi}.
The seventh and eighth parameters in Figure~\ref{fig:parameters-struct}
corresponds to $\beta$ described in Section~\ref{subsec:Symmetric-and-asymmetric}
and $\kappa_{\text{rrqr}}$ in Algorithm~\ref{alg:dense_factor},
respectively. We omit some additional parameters for simplicity.

\begin{figure}
\begin{lstlisting}[language={C++},basicstyle={\sffamily},tabsize=4]
struct hif::Params {
	double alpha_L;			// scalability oriented threshold for L (10)
	double alpha_U;			// scalability oriented threshold for U (10)
	double kappa;			// inverse norm threshold for L and U (3)
	double kappa_d;			// inverse norm threshold for D (3)
	double tau_L;			// drop tolerance for L (0.0001)
	double tau_U;			// drop tolerance for U (0.0001)
	double beta;			// safeguard for scaling in equilibration (1000)
	double kappa_rrqr;		// conditioning threshold for RRQR (eps^(-2/3))
};
\end{lstlisting}
\caption{\label{fig:parameters-struct}Core control parameters in HIFIR library
along with default values in parentheses.}
\end{figure}

The default values typically work well for PDE-based systems. They
are appropriate for both double- and single-precision computations,
since the thresholds are far greater than their corresponding machine
epsilons. We present some guidelines in tuning the first six parameters,
based on theoretical analysis and extensive experimentation:
\begin{itemize}
\item \textsf{alpha\_L}, \textsf{alpha\_U}: These control the scalability-oriented
dropping, and they are critical in achieving linear-time complexity
of HIFIR at each level. The recommended values are between $2$ and
$20$. While the default values $\alpha_{L}=\alpha_{U}=10$ are robust,
the user may reduce them to $\alpha_{L}=\alpha_{U}=3$ for PDE-based
systems (including saddle-point problems) for better efficiency;
\item \textsf{kappa, kappa\_d}: These control the inverse norms of the triangular
factors ($\boldsymbol{L}$ and $\boldsymbol{U}$) and diagonal factor
($\boldsymbol{D}$), respectively. Their recommended range is between
$3$ and $5$. The default values $\kappa=\kappa_{D}=3$ are robust,
but the user may increase them to $\kappa=\kappa_{D}=5$ for PDE-based
systems for better efficiency;
\item \textsf{tau\_L}, \textsf{tau\_U}: These control numerical droppings
and they are secondary compared to $\alpha_{L}$ and $\alpha_{U}$.
The default values $\tau_{L}=\tau_{U}=10^{-4}$ are robust; one may
increase them to $\tau_{L}=\tau_{U}=10^{-2}$ for PDE-based systems
for better efficiency. One could even set them to zero without losing
(near) linear time complexity, as long as $\alpha_{L}$ and $\alpha_{U}$
are within the recommended ranges. Note that large $\kappa$ and $\kappa_{D}$
may implicitly reduce the drop tolerances.
\end{itemize}

\subsection{Programming interfaces in C++}

The core components of HIFIR are implemented using generic object-oriented
programming with C++-11. All data structures and algorithms in HIFIR
are implemented as C++ template classes and functions.

\subsubsection{Basic data types}

HIFIR has two basic data types: sparse matrices and vectors. The sparse
matrices in HIF include CSC, CSR, and their partially or fully augmented
counterparts as described in Appendix~\ref{sec:Flexible-array-based}.
We will use CSR as the demonstration in the following. The class of
CSR matrix in HIF is hif::CSR<ValueType, IndexType>, of which the
two template arguments correspond to the value data type (e.g., \textsf{float},
\textsf{double}, \textsf{std::complex<double>}, etc.) and the index
data type (e.g., \textsf{int}, \textsf{long}, etc.), respectively.
For instance, a commonly used type is \textsf{hif::CSR<double, int>}.
A CSR instance can either wrap or own the data; the former mode allows
the user to create a ``view'' into a CSR matrix owned by another
software library, so it is more memory efficient and is preferred.
For example, to wrap external data for read-only access by HIF, one
can use the helper function
\begin{lstlisting}[language={C++},basicstyle={\sffamily},breaklines=true,tabsize=4]
	template<ValueType, IndexType> const hif::CSR<ValueType, IndexType> hif::wrap_const_csr(IndexType nrows, IndexType ncols, const IndexType *row_ptr, const IndexType *col_ind, const ValueType *vals);
\end{lstlisting}
where \textsf{row\_ptr}, \textsf{col\_ind}, and \textsf{vals} are
read-only. HIF also provides an interface to construct a CSR matrix
from scratch; we omit the details for simplicity. Note that the augmented
CSC and CSR data structures are used internally by HIF, so we omit
their interface definitions.

When using HIF as a preconditioner, the right-hand side and solution
vectors need to be passed in as arrays. HIF uses a generic interface
compatible with STL sequence containers such as std::vector. For memory
efficiency, HIF also provides a container \textsf{hif::Array<ValueType>}
for wrapping user-allocated arrays. The user can create mutable and
immutable instances by using one of the helper functions:
\begin{lstlisting}[language={C++},basicstyle={\sffamily},breaklines=true,tabsize=4]
	template<ValueType> hif::Array<ValueType> hif::wrap_array(std::size_t n, ValueType *vals);
	template<ValueType> const hif::Array<ValueType> hif::wrap_const_array(std::size_t n, const ValueType *vals);
\end{lstlisting}
If C++-20 is used, then \textsf{std::span<ValueType>} can be used
in place of \textsf{hif::Array<ValueType>}. Note that \textsf{hif::Array}
can also own the data; in this case, one can use a construct of \textsf{hif::Array}
similar to that of \textsf{std::vector}. 

\subsubsection{Interfaces for HIF preconditioners}

The algorithms for HIF are encapsulated in a class \textsf{hif::HIF<ValueType=double,
IndexType=int>}, of which the two template arguments are similar to
those of \textsf{hif::CSR} and specify the data types used by the
internal augmented CSC and CSR formats as well as the output format.
The class has two main interfaces: factorization and multilevel solve.
In the following, we will use \textsf{dHIF} as an alias for \textsf{hif::HIF<double,
int>} as a demonstration.

We first describe the member functions for computing the factorization,
which is provided by the template function
\begin{lstlisting}[language={C++},basicstyle={\sffamily},breaklines=true,tabsize=4,keepspaces=false]
	template<class MatType> void dHIF::factorize(const MatType &A, const hif::Params &params=hif::DEFAULT_PARAMS);
\end{lstlisting}
where \textsf{MatType} is typically \textsf{hif::CSR<ValueType, IndexType>}
or \textsf{hif::CSC<ValueType, IndexType>}. However, the interface
allows the use of different \textsf{ValueType} and \textsf{IndexType}
for \textsf{hif::HIF} and \textsf{MatType}. For example, the input
may be in double precision, and the preconditioner can be built in
single precision. However, if the input matrix is complex, then \textsf{ValueType}
of \textsf{hif::HIF} must also be complex and vice versa. The second
argument passes in the control parameters as described in Section~\ref{subsec:Control-parameters},
where the default is a static variable in \textsf{hif}.

The interface for the multilevel solve is provided by the template
function
\begin{lstlisting}[language={C++},basicstyle={\sffamily},breaklines=true,tabsize=4]
	template<class RhsType, class SolType> void dHIF::solve(const RhsType &b, SolType &x, bool trans=false, int rnk=0);
\end{lstlisting}
where \textsf{b} and \textsf{x} are the right-hand side and the solution
vector, respectively. The interfaces for \textsf{RhsType} and \textsf{SolType}
should be compatible with \textsf{const hif::Array<ValueType>} and
\textsf{hif::Array<ValueType>}, respectively, where their \textsf{ValueType}
may differ from each other and also differ from that of \textsf{hif::HIF}.
 The optional Boolean argument \textsf{trans} indicates whether to
apply the preconditioner itself or the (conjugate) transpose. The
argument \textsf{rnk} specifies the rank for the RRQR for the final
Schur complement; its default value $0$ indicates to use the rank
determined by \textsf{hif::HIF} from $\kappa_{\text{rrqr}}$. The
user can set \textsf{rnk} to $-1$ to use a larger rank determined
from $\kappa_{\text{rrqr}}=1/\epsilon_{\text{mach}}$ for the computation
of the null spaces. 

Besides the above main interfaces, \textsf{hif::HIF} also offers \textsf{hif::HIF::hifir}
for multilevel solve with iterative refinement and \textsf{hif::HIF::mmultiply}
for multilevel matrix-vector multiplication. Their interfaces are
as follows:

\begin{lstlisting}[language={C++},basicstyle={\sffamily},breaklines=true,tabsize=4]
	template<class MatType, class RhsType, class SolType> void dHIF::hifir(const MatType &A, const RhsType &b, SolType &x, int nirs, bool trans=false, int rnk=-1);
	template<class RhsType, class SolType> void dHIF::mmultiply(const RhsType &b, SolType &x, bool trans=false,  int rnk=-1);
\end{lstlisting}
These two functions are more useful for solving singular systems,
and hence their default values for \textsf{rnk} is $-1$ instead of
$0$. The argument \textsf{nirs} specifies the number of iterations
of iterative refinement. 

\subsection{\label{subsec:interface-for-Python-MATLAB}High-level interfaces
for MATLAB and Python}

It is often more productive for prototyping and academic research
to use a high-level programming language such as MATLAB, GNU Octave,
and Python. For this reason, we have developed \textsf{hifir4m} and
\textsf{hifir4py} to allow the users to access HIFIR from these languages. 

For the MATLAB interface, \textsf{hifir4m}, which also supports GNU
Octave, provides three key functions in its high-level programming
interface. The first function constructs a handle to a \textsf{HIF}
object for a matrix \textsf{A} from either \textsf{A} itself or an
optional ``sparsifier'' \textsf{S}, with the interface
\begin{lstlisting}[language=Matlab,basicstyle={\sffamily},breaklines=true,tabsize=4]
	function hif = hifCreate(A [, S, varargin])
\end{lstlisting}
where \textsf{A} and \textsf{S} can be MATLAB's built-in sparse format
matrices or MATLAB \textsf{struct} containing the three CSR fields,
and \textsf{varargin} specifies name-value pairs, such as \textsf{(...,
`alpha\_L', 5, `alpha\_U', 5, `mixed', true)}. The matrix \textsf{A}
and \textsf{S} can be single or double precision and can be real or
complex. The second function applies the preconditioner, with the
interface

\begin{lstlisting}[language=Matlab,basicstyle={\sffamily},breaklines=true,tabsize=4]
	function y = hifApply(hif, x [, op, rnk, nirs])
\end{lstlisting}
where the optional argument \textsf{op} can be one of `\textsf{S'},
`\textsf{SH'}, `\textsf{M'}, and `\textsf{MH'} for $\boldsymbol{M}^{g}\boldsymbol{x}$,
$\boldsymbol{M}^{gH}\boldsymbol{x}$, $\boldsymbol{M}\boldsymbol{x}$,
and $\boldsymbol{M}^{H}\boldsymbol{x}$, respectively, \textsf{rnk}
specifies the rank for truncated RRQR, and \textsf{nirs} specifies
the number of iterative refinements for $\boldsymbol{M}^{g}\boldsymbol{x}$
and $\boldsymbol{M}^{gH}\boldsymbol{x}$. To destroy the \textsf{HIF}
object, the user can call
\begin{lstlisting}[language=Matlab,basicstyle={\sffamily},breaklines=true,tabsize=4]
	delete(hif)
\end{lstlisting}
manually, or leave it to MATLAB to delete it automatically. In addition,
\textsf{hifir4m }offers two high-level drivers, \textsf{gmresHif}
and \textsf{pipitHifir}, for solving nonsingular and singular systems,
respectively. Their interfaces are similar to MATLAB's build-in \textsf{gmres}.

For the Python interface, \textsf{hifir4py} offers two sets of interfaces:
an intermediate-level interface consistent with the C++ version implemented
using Cython \cite{behnel2010cython}, and a high-level interface
consistent with the MATLAB version with support for SciPy sparse matrices
\cite{2020SciPy-NMeth}.

\section{\label{sec:Illustration-of-HIFIR}Illustrations of HIFIR for various
applications}

For HIFIR and its predecessor HILUCSI, we have reported extensive
comparisons with some prior state-of-the-art preconditioners \cite{chen2021hilucsi,chen2021robust,jiao2021approximate,jiao2021optimal}:
In \cite{chen2021hilucsi}, we compared HILUCSI with supernodal ILUTP
in SuperLU \cite{li2005overview,li2011supernodal} and multilevel
ILU in ILUPACK \cite{bollhofer2011ilupack} for indefinite systems;
in \cite{chen2021robust}, we compared HILUCSI with some customized
preconditioners (including \emph{pressure convection diffusion }(\emph{PCD})
\cite{silvester2001efficient,kay2002preconditioner}, \emph{least-squares
commutator} (\emph{LSC}) \cite{elman2006block}, and modified augmented
Lagrangian preconditioner \cite{benzi2006augmented,moulin2019augmented})
in Newton-GMRES method for solving the stationary incompressible Navier--Stokes
equations; in \cite{jiao2021optimal}, we compared HILUCSI with ILU(0)
\cite{saad2003iterative} and BoomerAMG in hypre \cite{hypre2021}
as building blocks for several block preconditioners for solving time-dependent
advection-diffusion equations with high-order finite element and finite
difference methods for spatial discretization and fully implicit multistage
Runge--Kutta schemes; in \cite{jiao2021approximate}, we compared
HIFIR with RIF-preconditioned LSMR \cite{benzi2003robust,fong2011lsmr}
and implicitly restarted Lanczos bidiagonalization \cite{baglama2005augmented}
for computing null-space vectors and for solving inconsistent singular
systems. In this section, we illustrate HIFIR for some of these applications,
using some of the aforementioned state-of-the-art techniques as points
of reference. We conducted our tests on a single node of a cluster
running CentOS 7.4 with dual 2.5 GHz 12-core Intel Xeon E5-2680v3
processors and 64 GB of RAM. We compiled HIFIR by GCC with optimization
options \textsf{-O3} and \textsf{-ffast-math}. For all the tests,
we used right-preconditioned GMRES with restart \cite[Section 9.3.2]{saad2003iterative},
and the dimension of KSP space is limited to $30$, i.e., GMRES($30$).

\subsection{Helmholtz equation}

As the first illustration, we solved the Helmholtz equation over $\Omega=\left[0,1\right]^{3}$
\begin{equation}
-\Delta u-k^{2}u=f\label{eq:helmholtz}
\end{equation}
with Dirichlet boundary conditions $u=u_{D}$ on $\partial\Omega$,
where $k>0$ is the wave number and $f$ is a source term. Such systems
are notoriously difficult to solve for have wave numbers \cite{ernst2012difficult}.
We discretized the equation using the Galerkin finite element methods
(FEM) with quadratic tetrahedral (aka $P_{2}$) elements, which we
implemented using FEniCS v2019.1.0 \cite{logg2010dolfin,alnaes2015fenics}.
We constructed the right-hand side and the boundary conditions using
the method of manufactured solutions with the exact solution $u(x,y,z)=\cos(\pi x)\sin(\pi y)\sin(\pi z)$.
Given a consistent numerical discretization method for (\ref{eq:helmholtz}),
we arrive at the following system of equations
\begin{equation}
\left(\boldsymbol{K}-k^{2}\boldsymbol{M}\right)\boldsymbol{u}_{h}=\boldsymbol{f}_{h},\label{eq:helm-discretize}
\end{equation}
where $\boldsymbol{K}$ and $\boldsymbol{M}$ are the stiffness and
mass matrices, respectively, and $\boldsymbol{u}_{h}$ and $\boldsymbol{f}_{h}$
are vectors containing nodal values of $u$ and $f$, respectively.
Note that $\boldsymbol{K}$ and $\boldsymbol{M}$ are (symmetric)
positive definite, but the coefficient matrix $\boldsymbol{K}-k^{2}\boldsymbol{M}$
can be positive definite, indefinite, or even singular.  We chose
three different wave numbers, $k=1$, $5$, and $10$. The first two
resulted in positive-definite systems, and the third resulted in indefinite
systems.

To assess the effectiveness of HIF, we assembled three systems using
meshes with 93,750, 257,250, and 795,906 quadratic elements, respectively.
 We solved the systems using $\text{rtol}=10^{-6}$ GMRES(30) and
$\tau=10^{-2}$, $\kappa=5$, and $\alpha=3$ for double-precision
HIF. For each of these systems, HIF ended up producing a four-level
factorization with nnz ratios (i.e., the number of nonzeros in the
preconditioner versus that in the input matrix) of about 2.7. As points
of reference, we also solved the systems using GMRES(30) with the
supernodal ILUTP in SuperLU v5.2.2 \cite{li2011supernodal,li2005overview}
and ILUTP in WSMP v20.12 \cite{gupta2001wsmp,gupta2021wsmp} as the
right preconditioners with $\text{droptol}=10^{-2}$. As can be seen
from Table~\ref{tab:res-helm}, HIF was the fastest for seven out
of nine cases, and WSMP was the fastest for two small cases. SuperLU
was substantially slower than both HIF and WSMP. Nevertheless, SuperLU
was relatively insensitive to the wave numbers in terms of factorization
and total times, as was HIF. In contrast, WSMP had a significant increase
in the solve times for $k=10$, indicating that ILUTP in WSMP is less
optimized for indefinite systems than for positive-definite systems.
In addition, the computational cost of WSMP increased at a faster
rate than HIF as the problem size increased. As a result, HIF was
about a factor of 2.2 faster than WSMP for $k=10$ on the finest mesh.
The better performance of HIF was mostly due to its scalability-oriented
dropping. It is worth noting that WSMP supports parallel ILUTP and
all of our comparisons were conducted in serial. Parallel ILUTP in
WSMP may potentially outperform HIF on a multicore computer, but parallel
ILUTP is less robust than serial ILUTP in WSMP and caused GMRES to
stagnate for larger systems in our tests.

\begin{table}
\caption{\label{tab:res-helm}Results for systems arising from the Helmholtz
equation discretized by $P_{2}$ FEM, solved using GMRES(30) with
HIF compared to with supernodal ILUTP in SuperLU and ILUTP in WSMP
as the right preconditioner. The condition numbers were estimated
using \textsf{condest} function in MATLAB. Times are in seconds. Leaders
are in boldface.}

\centering{}\setlength\tabcolsep{3pt}%
\begin{tabular}{ccccc|ccc|ccc|ccc}
\hline 
\multirow{2}{*}{mesh} & \multirow{2}{*}{$n$} & \multirow{2}{*}{nnz} & \multirow{2}{*}{$k$} & \multirow{2}{*}{cond.} & \multicolumn{3}{c|}{HIF} & \multicolumn{3}{c|}{supernodal ILUTP } & \multicolumn{3}{c}{WSMP/ILUTP}\tabularnewline
\cline{6-14} \cline{7-14} \cline{8-14} \cline{9-14} \cline{10-14} \cline{11-14} \cline{12-14} \cline{13-14} \cline{14-14} 
 &  &  &  &  & $\begin{array}{c}
\text{fac.}\\
\text{time}
\end{array}$ & $\begin{array}{c}
\text{tot.}\\
\text{time}
\end{array}$ & iter. & $\begin{array}{c}
\text{fac.}\\
\text{time}
\end{array}$ & $\begin{array}{c}
\text{tot.}\\
\text{time}
\end{array}$ & iter. & $\begin{array}{c}
\text{fac.}\\
\text{time}
\end{array}$ & $\begin{array}{c}
\text{tot.}\\
\text{time}
\end{array}$ & iter.\tabularnewline
\hline 
\hline 
\multirow{3}{*}{coarse} & \multirow{3}{*}{132,651} & \multirow{3}{*}{3,195,529} & 1 & 8.8e4 & 3.15 & 3.46 & \textbf{8} & 20.2 & 21.1 & 11 & \textbf{2.25} & \textbf{2.91} & 23\tabularnewline
 &  &  & 5 & 6.4e5 & 3.1 & 3.47 & \textbf{8} & 20.9 & 21.6 & 12 & \textbf{2.26} & \textbf{3.37} & 36\tabularnewline
 &  &  & 10 & 6.5e5 & 3.03 & \textbf{3.81} & \textbf{22} & 19.8 & 21.8 & 25 & \textbf{2.27} & 5.42 & 103\tabularnewline
\hline 
\multirow{3}{*}{medium} & \multirow{3}{*}{357,911} & \multirow{3}{*}{9,070,749} & 1 & 2.4e5 & 8.82 & \textbf{9.62} & \textbf{8} & 115 & 117 & 12 & \textbf{7.85} & 10.8 & 30\tabularnewline
 &  &  & 5 & 1.8e6 & 9.02 & \textbf{9.89} & \textbf{9} & 117 & 120 & 13 & \textbf{7.86} & 11.9 & 34\tabularnewline
 &  &  & 10 & 1.8e6 & 8.9 & \textbf{10.7} & \textbf{16} & 116 & 120 & 24 & \textbf{7.9} & 25.1 & 192\tabularnewline
\hline 
\multirow{3}{*}{fine} & \multirow{3}{*}{1,092,727} & \multirow{3}{*}{28,814,525} & 1 & 7.4e5 & \textbf{29.8} & \textbf{33.8} & \textbf{11} & 903 & 941 & 15 & 35.5 & 46.1 & 32\tabularnewline
 &  &  & 5 & 5.4e6 & \textbf{30.1} & \textbf{34.1} & \textbf{11} & 918 & 944 & 16 & 35.4 & 47.9 & 37\tabularnewline
 &  &  & 10 & 5.5e6 & \textbf{29.8} & \textbf{35.5} & \textbf{16} & 903 & 954 & 20 & 35.1 & 77.4 & 141\tabularnewline
\hline 
\end{tabular}
\end{table}

\subsection{Linear elasticity with pure traction boundary conditions}

In this illustration, we considered the linear-elasticity model of
a solid body $\Omega\subset\mathbb{R}^{3}$,
\begin{equation}
-\boldsymbol{\nabla}\cdot\boldsymbol{\sigma}=\boldsymbol{f},\label{eq:le}
\end{equation}
where $\boldsymbol{\sigma}$ and $\boldsymbol{f}$ are Cauchy stress
tensor and body force per unit volume, respectively. For isotropic
material, $\boldsymbol{\sigma}=\lambda(\boldsymbol{\nabla}\cdot\boldsymbol{u})\boldsymbol{I}+\mu(\boldsymbol{\nabla}\boldsymbol{u}+(\boldsymbol{\nabla}\boldsymbol{u})^{T})$,
where $\boldsymbol{u}$ is the displacement, and $\lambda$ and $\mu$
are Lam\'e's parameters. We applied pure traction (aka Neumann) boundary
conditions to (\ref{eq:le}), i.e.,
\begin{equation}
\boldsymbol{\sigma}\cdot\boldsymbol{n}=\boldsymbol{t}\qquad\text{on }\partial\Omega,\label{eq:pure-traction}
\end{equation}
where $\boldsymbol{n}$ and $\boldsymbol{t}$ are unit outward surface
normal and surface traction, respectively. Discretizing (\ref{eq:le})
and (\ref{eq:pure-traction}) using a Galerkin FEM leads to a singular
system
\begin{equation}
\boldsymbol{K}\boldsymbol{u}_{h}=\boldsymbol{f}_{h},\label{eq:inconsist-le-sys}
\end{equation}
where the stiffness matrix $\boldsymbol{K}\in\mathbb{R}^{n\times n}$
is symmetric positive semidefinite (SPSD), and $\boldsymbol{u}_{h}$
and $\boldsymbol{f}_{h}$ contain nodal values of $\boldsymbol{u}$
and $f$, respectively. The continuum equation is invariant to translation
and rotation, so the null space of (\ref{eq:inconsist-le-sys}) is
six dimensional. As in \cite{kuchta2019singular}, we constructed
the solid domain $\Omega$ by first rotating a box $\left[\nicefrac{-1}{4},\nicefrac{1}{4}\right]\times\left[\nicefrac{-1}{2},\nicefrac{1}{2}\right]\times\left[\nicefrac{-1}{8},\nicefrac{1}{8}\right]$
around $x$-,\textbf{ }$y$-, and $z$-axes by $\nicefrac{\pi}{2}$,
$\nicefrac{\pi}{4}$, $\nicefrac{\pi}{5}$ in that order and then
by translating it by $\left[0.1,0.2,0.3\right]{}^{T}$. We discretized
the domain using three linear tetrahedral meshes with 24,576, 196,608,
and 2,058,000 $P_{1}$ elements, respectively. We computed the body
force $\boldsymbol{f}$ and surface traction $\boldsymbol{h}$ using
the manufactured solution $\boldsymbol{u}(x,y,z)=\frac{1}{4}(\sin(\frac{\pi}{4}x),z^{3},-y)$.
Due to discretization errors, the resulting linear system (\ref{eq:inconsist-le-sys})
may be inconsistent, which poses significant challenges for computing
a stable solution for $\boldsymbol{u}_{h}$. Note that the nullspace
of $\boldsymbol{K}$ corresponds the rigid-body motion, of which an
orthonormal basis may be computed analytically if the mesh and the
finite-element basis functions are known; see e.g., \cite{kuchta2019singular}.
However, an algebraic solver should be able to compute the nullspace
from the matrix directly, for example, when the solver does not have
access to the mesh or the basis functions.

To solve the inconsistent system, we use the solver \emph{PIPIT} (or
\emph{PseudoInverse solver via Preconditioned ITerations}) \cite{jiao2021approximate},
which is based on HIFIR and is composed of three steps:
\begin{enumerate}
\item Compute orthonormal basis $\boldsymbol{V}\in\mathbb{R}^{n\times6}$
of $\mathcal{N}(\boldsymbol{K}^{T})\equiv\mathcal{N}(\boldsymbol{K})$
for $\boldsymbol{K}$ in (\ref{eq:inconsist-le-sys}) using HIFIR-preconditioned
FGMRES.
\item Find a least-squares solution $\boldsymbol{u}_{\text{LS}}$ of the
consistent system $\boldsymbol{K}\boldsymbol{u}_{\text{LS}}=(\boldsymbol{I}-\boldsymbol{V}\boldsymbol{V}^{T})\boldsymbol{f}_{h}$
using HIF-preconditioned GMRES.
\item Obtain the pseudoinverse solution $\boldsymbol{u}_{\text{PI}}$ via
orthogonal projection, i.e., $\boldsymbol{u}_{\text{PI}}=(\boldsymbol{I}-\boldsymbol{V}\boldsymbol{V}^{T})\boldsymbol{u}_{\text{LS}}$.
\end{enumerate}
All three steps reuse the same incomplete factorization. We refer
readers to \cite{jiao2021approximate} for more detail. If the nullspace
is known \emph{a priori}, such as in \cite{kuchta2019singular}, then
Step 1 can be omitted. When computing the null-space vectors, we used
FGMRES(30) with Householder QR for the Arnoldi process and terminated
when the Hessenberg matrix became ill-conditioned. We used the default
parameters (see Section~\ref{subsec:Control-parameters}) for HIFIR
except that $\text{rtol}$ was set to $10^{-10}$ in step 2. For these
meshes, HIF ended up having three, four, and five levels with nnz
ratios 8.84, 10.4, and 10.7, respectively. Table~\ref{tab:res-le}
summarizes the timing results of PIPIT for these singular systems.
As a point of reference, we used LSMR \cite{fong2011lsmr} preconditioned
with RIF \cite{benzi2003robust}. Since LSMR did not converge in our
test without a preconditioner and a preconditioned LSMR can only compute
a least-squares solution, we compare the computational cost of PIPIT
for computing a least-squares solution (i.e., the factorization cost
of HIF and step 2 of PIPIT) with RIF+LSMR. It can be seen that PIPIT
was more than an order of magnitude more efficient than RIF+LSMR for
the small and medium meshes. For the large mesh, RIF+LSMR failed to
converge to the desired precision after 10,000 iterations. In addition,
Table~\ref{tab:res-le} also shows the runtimes for step 1 of PIPIT,
which were comparable to solving for $\boldsymbol{u}_{\text{LS}}$
since the same HIF preconditioner was reused. In addition, we report
the accuracy of the first ($\boldsymbol{v}_{1}$) and sixth ($\boldsymbol{v}_{6}$)
null-space vectors as examples, both of which converted to (near)
machine precision. Note that for RIF+LSMR to compute $\boldsymbol{u}_{\text{PI}}$
instead of $\boldsymbol{u}_{\text{LS}}$, it would also need to compute
the null-space vectors, presumably using preconditioned LSMR \cite{fong2011lsmr},
which could not converge at least for the largest system.

\begin{table}
\caption{\label{tab:res-le}Results for singular systems arising from linear
elasticity with pure traction boundary conditions using $P_{1}$ finite
elements. We compare the costs of step 2 of PIPIT with RIF+LSMR, where
leaders are in boldface and `$-$' indicated non-convergence after
10,000 iterations. The last three columns show the cost and accuracy
of null-space computations in step 1 of PIPIT. Times are in seconds.}

\centering{}\setlength\tabcolsep{2pt}%
\begin{tabular}{ccc||ccc|ccc||ccc}
\hline 
\multirow{2}{*}{mesh} & \multirow{2}{*}{$n$} & \multirow{2}{*}{nnz} & \multicolumn{3}{c|}{PIPIT for $\boldsymbol{u}_{\text{LS}}$} & \multicolumn{3}{c||}{RIF+LSMR for $\boldsymbol{u}_{\text{LS}}$} & \multicolumn{3}{c}{null-space computation}\tabularnewline
\cline{4-12} \cline{5-12} \cline{6-12} \cline{7-12} \cline{8-12} \cline{9-12} \cline{10-12} \cline{11-12} \cline{12-12} 
 &  &  & fac. time & tot. time & iter. & fac. time & tot. time & iter. & time & $\nicefrac{\left\Vert \boldsymbol{K}\boldsymbol{v}_{1}\right\Vert }{\left\Vert \boldsymbol{K}\right\Vert }$ & $\nicefrac{\left\Vert \boldsymbol{K}\boldsymbol{v}_{6}\right\Vert }{\left\Vert \boldsymbol{K}\right\Vert }$\tabularnewline
\hline 
\hline 
\multirow{1}{*}{coarse} & 15,147 & \multirow{1}{*}{610,929} & \multirow{1}{*}{\textbf{3.5}} & \multirow{1}{*}{\textbf{3.82}} & \multirow{1}{*}{\textbf{15}} & \multirow{1}{*}{12.5} & \multirow{1}{*}{30.5} & \multirow{1}{*}{2,588} & 3.01 & 7e-16 & 3e-14\tabularnewline
\multirow{1}{*}{medium} & 109,395 & \multirow{1}{*}{4,652,505} & \multirow{1}{*}{\textbf{39.4}} & \multirow{1}{*}{\textbf{44}} & \multirow{1}{*}{\textbf{35}} & \multirow{1}{*}{308} & \multirow{1}{*}{741} & \multirow{1}{*}{8,673} & 126 & 3e-15 & 4e-16\tabularnewline
\multirow{1}{*}{fine} & 1,081,188 & \multirow{1}{*}{47,392,074} & \multirow{1}{*}{\textbf{511}} & \multirow{1}{*}{\textbf{700}} & \multirow{1}{*}{\textbf{152}} & \multirow{1}{*}{$-$} & \multirow{1}{*}{$-$} & \multirow{1}{*}{$-$} & 3.1e3 & 4e-15 & 6e-16\tabularnewline
\hline 
\end{tabular}
\end{table}

\subsection{Stationary incompressible Navier--Stokes equations}

As a third illustration, we solve the incompressible Navier--Stokes
(INS) equations for modeling fluids. We consider the INS equations
with normalized density on a domain $\Omega\subset\mathbb{R}^{3}$,
which read
\begin{align}
\boldsymbol{u}_{t}-\nu\Delta\boldsymbol{u}+\boldsymbol{u}\cdot\boldsymbol{\nabla u}+\boldsymbol{\nabla}p & =\boldsymbol{g},\label{eq:momentum}\\
\boldsymbol{\nabla}\cdot\boldsymbol{u} & =0,\label{eq:div-free}
\end{align}
where $\boldsymbol{u}$ and $p$ are velocities and pressure, respectively,
$\nu$ is the kinetic viscosity, and the subscript $t$ denotes the
temporal derivative. After dropping off the temporal-derivative term,
we arrive at the stationary INS, for which the momentum equation (\ref{eq:momentum})
becomes
\begin{equation}
-\nu\Delta\boldsymbol{u}+\boldsymbol{u}\cdot\boldsymbol{\nabla u}+\boldsymbol{\nabla}p=\boldsymbol{g}.\label{eq:steady-ins}
\end{equation}
When integrating (\ref{eq:momentum}) in time, (\ref{eq:steady-ins})
is equivalent to having the time step equal to infinity. We discretize
the INS in space using the $P_{2}$-$P_{1}$ Taylor--Hood (TH) elements
\cite{taylor1973numerical}, which leads to a system of nonlinear
systems of equations. We then solve it using Newton-GMRES, an \emph{inexact
Newton method} \cite{dembo1982inexact} that uses the preconditioned
GMRES in the inner iterations of Newton's method. At iteration $k$,
Newton's method solves a linear problem
\begin{equation}
\boldsymbol{J}_{k}\begin{bmatrix}\delta\boldsymbol{u}_{k}\\
\delta\boldsymbol{p}_{k}
\end{bmatrix}\approx-\begin{bmatrix}\boldsymbol{f}_{k}\\
\boldsymbol{y}_{k}
\end{bmatrix}\qquad\text{with}\qquad\boldsymbol{J}_{k}=\begin{bmatrix}\boldsymbol{K}+\boldsymbol{C}_{k}+\boldsymbol{W}_{k} & \boldsymbol{E}^{T}\\
\boldsymbol{E} & \boldsymbol{0}
\end{bmatrix},\label{eq:ins-nt}
\end{equation}
where $\boldsymbol{K}$, $\boldsymbol{E}$, $\boldsymbol{C}_{k}$,
and $\boldsymbol{W}_{k}$ correspond to $\nu\Delta\delta\boldsymbol{u}$,
$\boldsymbol{\nabla}\cdot\delta\boldsymbol{u}$, $\boldsymbol{u}_{k}\cdot\boldsymbol{\nabla}\delta\boldsymbol{u}$,
and $\delta\boldsymbol{u}\cdot\boldsymbol{\nabla}\boldsymbol{u}_{k}$,
respectively, and $\delta\boldsymbol{u}_{k}$ and $\delta\boldsymbol{p}_{k}$
denote increments in velocity and pressure, respectively.  Since
$\boldsymbol{J}_{k}$ is fairly dense in Newton's method, we constructed
a ``sparsifer'' $\boldsymbol{B}_{k}$ by omitting $\boldsymbol{W}_{k}$
when building the preconditioner. For robustness, we start with Picard
iterations, which uses $\boldsymbol{B}_{k}$ in place of $\boldsymbol{J}_{k}$
at each Newton's step, and then switch to Newton's iterations once
the solution is sufficiently accurate. For the initial guess $\left[\boldsymbol{u}_{0},\boldsymbol{p}_{0}\right]^{T}$,
we obtained it by solving the corresponding Stokes equation, i.e.,
$-\nu\Delta\boldsymbol{u}+\boldsymbol{\nabla}p=\boldsymbol{g}$.
For more details, see \cite{chen2021robust}.

\begin{figure}
\begin{centering}
\includegraphics[width=0.7\columnwidth]{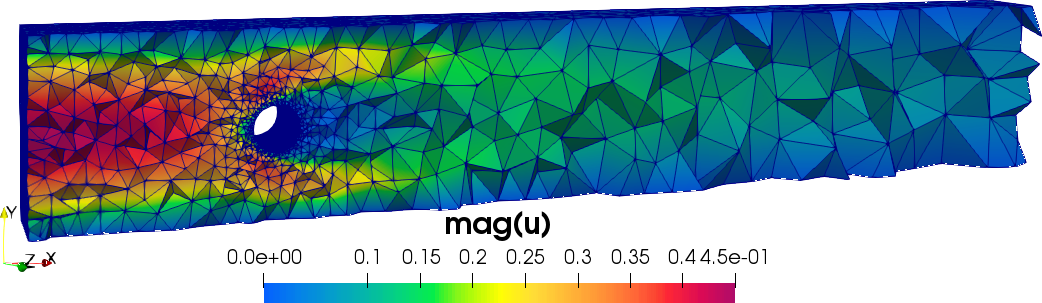}
\par\end{centering}
\caption{\label{fig:3D-flow-over-cylinder}Cut-off view of the speed for the
3D flow-over-cylinder problem with a coarse sample mesh.}
\end{figure}

As a demonstration, we solved the 3D flow-over-cylinder benchmark
problem \cite{schafer1996benchmark}, where a cylinder of diameter
$D=0.1$ is placed in a channel of length $2.5$ with square cross-sections
of height $H=0.41$. The inflow boundary condition is imposed on the
left side of a box, which reads $\boldsymbol{u}_{\text{in}}=(U(y,z),0,0)$,
where $U(y,z)=16\times0.45\thinspace yz(H-y)(H-z)/H^{4}$. A ``do-nothing''
velocity is imposed for the outflow (right face) along with a zero
pressure. The rest boundaries are no-slip walls. The kinetic viscosity
was $\nu=1\times10^{-3}$, so the Reynolds number is $\text{Re}=4\times0.45\thinspace D/(9\thinspace\nu)=20$.
We generated three sets of unstructured tetrahedral meshes using Gmsh
\cite{geuzaine2009gmsh}, with 71,031, 268,814, and 930,248 elements,
respectively. Figure~\ref{fig:3D-flow-over-cylinder} shows a coarse
sample mesh along with the computed speed.
 We assembled the matrices $\boldsymbol{J}_{k}$ and $\boldsymbol{B}_{k}$
using our in-house FEM code. Table~\ref{tab:res-ins} shows the performance
of HIF+GMRES(30) for a representative Newton step, with $\text{rtol}=10^{-6}$
in GMRES and $\tau=10^{-2}$, $\kappa=5$, and $\alpha=3$ for HIF.
As a point of reference, we solved the same systems using ILUPACK
\cite{bollhofer2011ilupack} with $\text{droptol}=10^{-2}$ and other
default parameters. As we can see, HIF-preconditioned GMRES performed
consistently well, even for the finest mesh. Its overall performance
was faster than ILUPACK by a factor about six and 16 for the coarse
and intermediate-level meshes, although ILUPACK used fewer GMRES iterations.
More importantly, ILUPACK ran out of the 64 GB main memory for the
largest problem.  In \cite{chen2021robust}, we also compared HIF
with the direct solve MUMPS \cite{amestoy2000mumps}, which also ran
out of memory for the largest case as ILUPACK, and with ILU($1$)
and ILU($2$), for which GMRES failed to converge. 

\begin{table}
\caption{\label{tab:res-ins}Results for stationary INS with TH elements in
one representative Newton's step solved by GMRES(3) preconditioned
by HIF on the sparsifier $\boldsymbol{B}_{k}$ in comparison with
ILUPACK. `$\times$' indicates out of memory. Times are in seconds.
Leaders are in boldface.}

\centering{}\setlength\tabcolsep{3pt}%
\begin{tabular}{cccc|cccc|cccc}
\hline 
\multirow{2}{*}{mesh} & \multirow{2}{*}{$n$} & \multirow{2}{*}{$\text{nnz}(\boldsymbol{J}_{k})$} & \multirow{2}{*}{$\text{nnz}(\boldsymbol{B}_{k})$} & \multicolumn{4}{c|}{HIF+GMRES(30)} & \multicolumn{4}{c}{ILUPACK+GMRES(30)}\tabularnewline
\cline{5-12} \cline{6-12} \cline{7-12} \cline{8-12} \cline{9-12} \cline{10-12} \cline{11-12} \cline{12-12} 
 &  &  &  & nnz rat. & fac. time & tot. time & iter. & nnz rat. & fac. time & tot. time & iter.\tabularnewline
\hline 
\hline 
\multirow{1}{*}{coarse} & \multirow{1}{*}{262,912} & 21,870,739 & \multirow{1}{*}{$\text{9,902,533}$} & \multirow{1}{*}{3.47} & \multirow{1}{*}{\textbf{20.3}} & \multirow{1}{*}{\textbf{22.6}} & \multirow{1}{*}{17} & \multirow{1}{*}{8.34} & \multirow{1}{*}{130} & \multirow{1}{*}{133} & \multirow{1}{*}{\textbf{11}}\tabularnewline
\multirow{1}{*}{medium} & \multirow{1}{*}{1,086,263} & 98,205,997 & \multirow{1}{*}{43,686,979} & \multirow{1}{*}{3.4} & \multirow{1}{*}{\textbf{106}} & \multirow{1}{*}{\textbf{124}} & \multirow{1}{*}{33} & \multirow{1}{*}{13.4} & \multirow{1}{*}{1.92e3} & \multirow{1}{*}{1.94e3} & \multirow{1}{*}{\textbf{14}}\tabularnewline
\multirow{1}{*}{fine} & \multirow{1}{*}{3,738,327} & 343,357,455 & \multirow{1}{*}{152,438,721} & \multirow{1}{*}{3.47} & \multirow{1}{*}{\textbf{442}} & \multirow{1}{*}{\textbf{586}} & \multirow{1}{*}{\textbf{68}} & \multirow{1}{*}{$\times$} & \multirow{1}{*}{$\times$} & \multirow{1}{*}{$\times$} & \multirow{1}{*}{$\times$}\tabularnewline
\hline 
\end{tabular}
\end{table}

\subsection{Advection-diffusion equation with fully implicit Runge--Kutta schemes}

As our final illustration, we consider the time-dependent advection-diffusion
(AD) equation for $u\colon\Omega\times[0,T]\rightarrow\mathbb{R}$,
\begin{equation}
u_{t}-\mu\Delta u+\boldsymbol{v}\cdot\boldsymbol{\nabla}u=f,\label{eq:time-ad-eqn}
\end{equation}
where $\mu\ge0$ is the diffusion coefficient, $\boldsymbol{v}$ is
a divergence-free velocity field, and $f$ is some source term. As
a demonstration, we chose $\Omega=\left[0,1\right]^{3}$ and computed
$f$ and the Dirichlet boundary conditions from the manufactured solution
$u(x,y,z,t)=\sin\left(1.5\pi t\right)\sin\left(\pi x\right)\sin\left(\pi y\right)\sin\left(\pi z\right)$.
We used the fourth and sixth-order finite difference method (FDM)
to discretize (\ref{eq:time-ad-eqn}), leading to the semi-discretization
form
\begin{equation}
\boldsymbol{u}_{t}\left(t\right)=\boldsymbol{J}\boldsymbol{u}\left(t\right)+\boldsymbol{f}\left(t\right),\label{eq:ad-odes}
\end{equation}
where $\boldsymbol{J}\in\mathbb{R}^{n\times n}$ denotes the Jacobian
matrix corresponding to the operator $(\mu\Delta-\boldsymbol{v}\cdot\boldsymbol{\nabla})$.
Eq.~(\ref{eq:ad-odes}) is a system of stiff ordinary differential
equations (ODEs) and can be integrated using high-order fully implicit
Runge--Kutta (FIRK) schemes, such as $s$-stage Gauss--Legendre
schemes, which are $2s$-order accurate and are unconditionally stable.
A Runge--Kutta scheme can be expressed by the Butcher tableau $\begin{array}{c|c}
\boldsymbol{c} & \boldsymbol{A}\\
\hline  & \boldsymbol{b}^{T}
\end{array}$, where $\boldsymbol{A}\in\mathbb{R}^{s\times s}$, $\boldsymbol{c}\in\mathbb{R}^{s}$,
and $\boldsymbol{b}\in\mathbb{R}^{s}$ \cite{butcher2008numerical}.
Given a time step $\delta t$, Eq.~(\ref{eq:ad-odes}) leads to an
$sn\times sn$ linear system
\begin{equation}
\boldsymbol{\mathcal{A}}\mathcal{K}=\mathcal{B}\qquad\text{with}\qquad\boldsymbol{\mathcal{A}}=\boldsymbol{I}_{sn}-\delta t\boldsymbol{A}\otimes\boldsymbol{J},\label{eq:rk-A}
\end{equation}
where $\boldsymbol{I}_{m}$ denotes the $m$-dimension identity matrix,
$\otimes$ denotes the Kronecker-product operator, $\boldsymbol{\mathcal{A}}\in\mathbb{R}^{sn\times sn}$,
$\mathcal{B}\in\mathbb{R}^{sn}$, and $\mathcal{K}\in\mathbb{R}^{sn}$.
We refer readers to \cite{jiao2021optimal} for more details. 

To assess HIF for complex-valued matrices, we use an optimal preconditioner,
called \emph{block CSD} or \emph{BCSD}, developed in \cite{jiao2021optimal}.
BCSD is based on the complex Schur decomposition (CSD) of the Butcher
matrix $\boldsymbol{A}$ (i.e., $\boldsymbol{A}=\boldsymbol{Q}\boldsymbol{U}\boldsymbol{Q}^{H}$,
where $\boldsymbol{Q}\in\mathbb{C}^{s\times s}$ is unitary and $\boldsymbol{U}\in\mathbb{C}^{s\times s}$
is upper triangular), and it reads
\begin{equation}
\boldsymbol{\mathcal{M}}=\left(\boldsymbol{Q}\otimes\boldsymbol{I}_{n}\right)\boldsymbol{\mathcal{U}}\left(\boldsymbol{Q}^{H}\otimes\boldsymbol{I}_{n}\right)\qquad\text{with}\qquad\boldsymbol{\mathcal{U}}=\boldsymbol{I}_{sn}-\delta t\boldsymbol{U}\otimes\boldsymbol{J}.\label{eq:bcsd}
\end{equation}
$\boldsymbol{\mathcal{M}}$ is optimal in that $\boldsymbol{\mathcal{A}}\boldsymbol{\mathcal{M}}^{-1}=\boldsymbol{I}_{sn}$,
where $\boldsymbol{\mathcal{M}}^{-1}=(\boldsymbol{Q}\otimes\boldsymbol{I}_{n})\boldsymbol{\mathcal{U}}^{-1}(\boldsymbol{Q}^{H}\otimes\boldsymbol{I}_{n})$;
see \cite{jiao2021optimal} for more details. Note that $\boldsymbol{\mathcal{U}}$
in (\ref{eq:bcsd}) is block upper triangular, so computing $\boldsymbol{\mathcal{U}}^{-1}\mathcal{X}$
for a vector $\mathcal{X}$ only requires (approximately) factorizing
the diagonal blocks of $\boldsymbol{\mathcal{U}}$, i.e., $\boldsymbol{M}_{i}=\boldsymbol{I}_{n}-\delta t\thinspace\lambda_{i}\boldsymbol{J}\in\mathbb{C}^{n\times n}$
for $i=1,2,\ldots s,$ where $\lambda_{i}\in\mathbb{C}$ are the diagonal
entries of $\boldsymbol{U}$, which are the eigenvalues of $\boldsymbol{A}$.
For Gauss--Legendre schemes, there are $\left\lfloor \nicefrac{s}{2}\right\rfloor $
distinct conjugate pairs of complex eigenvalues of $\boldsymbol{A}$,
and given $\boldsymbol{M}_{i}\in\mathbb{C}^{n\times n}$, $\boldsymbol{y}=\overline{\boldsymbol{M}_{i}}^{-1}\boldsymbol{x}\implies\overline{\boldsymbol{y}}=\overline{\boldsymbol{M}_{i}^{-1}\boldsymbol{x}}$.
Hence, we only need to factorize $\left\lceil \nicefrac{s}{2}\right\rceil $
distinct complex-valued diagonal blocks for BCSD. We constructed the
blocks in the Jacobian matrices $\boldsymbol{J}$ in (\ref{eq:ad-odes})
using our in-house high-order FDM code on equidistant structured grids
with grid sizes $h=\nicefrac{1}{32},\nicefrac{1}{64},\nicefrac{1}{128}$
and a velocity field $\boldsymbol{v}=[1,1,1]^{T}$. We factorized
$\boldsymbol{M}_{i}$ using double-precision complex HIF with parameters
$\tau=10^{-3}$, $\kappa=5$, and $\alpha=5$. With these parameters,
HIF ended up producing three to four levels. Table~\ref{tab:res-ad-rk}
shows the performance of HIF-based BCSD preconditioner for the two-
and four-stage Gauss--Legendre schemes for the first time step with
$\delta t=\nicefrac{1}{8}$. As a point of reference, Table~\ref{tab:res-ad-rk}
also reports the results using MATLAB's built-in (precompiled) \textsf{ilu}
function without fills (aka ILU(0)) with preprocessing steps (including
equilibration \cite{duff2001algorithms} and fill-reduction reordering
\cite{chan1980linear}). We chose ILU(0) as a baseline since it is
commonly used in the literature for the diagonal blocks \cite{kanevsky2007application,pazner2017stage}.
It can be seen HIF was significantly faster than ILU(0) overall, although
ILU(0) was more efficient in terms of factorization cost. In addition,
BCSD with ILU(0) failed to converge for the largest case. In \cite{jiao2021optimal},
we also compared HIF with BoomerAMG in hypre v2.21.0 \cite{hypre2021}
for some real-valued block preconditioners, and HIF was about a factor
of eight faster than BoomerAMG for larger problems. We do not compare
BoomerAMG for BCSD, since hypre does not yet support complex arithmetic.

\begin{table}
\caption{\label{tab:res-ad-rk}Results for solving the advection-diffusion
equation with fourth- and sixth-order FDM and two- and four-stage
Gauss--Legendre schemes correspondingly, using GMRES(30) right-preconditioned
by HIF-based and ILU(0)-based BCSD. Times are in seconds. `$-$' indicates
that GMRES failed to converge after 500 iterations. Leaders are in
boldface.}

\centering{}%
\begin{tabular}{ccc||cccc|cccc}
\hline 
\multirow{2}{*}{$h$} & \multirow{2}{*}{$n$} & \multirow{2}{*}{nnz} & \multicolumn{4}{c|}{BCSD with HIF} & \multicolumn{4}{c}{BCSD with ILU(0)+pre}\tabularnewline
\cline{4-11} \cline{5-11} \cline{6-11} \cline{7-11} \cline{8-11} \cline{9-11} \cline{10-11} \cline{11-11} 
 &  &  & nnz ratio & fac. time & tot. time & iter. & nnz ratio & fac. time & tot. time & iter.\tabularnewline
\hline 
\hline 
\multicolumn{11}{c}{Two-stage Gauss--Legendre scheme with fourth-order finite differences}\tabularnewline
\hline 
$\nicefrac{1}{32}$ & 29,791 & 381,517 & 5.47 & 1.08 & 1.52 & \textbf{17} & 1.0 & \textbf{0.1} & \textbf{1.31} & 66\tabularnewline
$\nicefrac{1}{64}$ & 250,047 & 3,226,797 & 5.27 & 11.0 & \textbf{19.5} & \textbf{34} & 1.0 & \textbf{0.89} & 23.8 & 161\tabularnewline
$\nicefrac{1}{128}$ & 2,048,383 & 26,532,205 & 5.27 & 104 & \textbf{262} & \textbf{77} & 1.0 & \textbf{8.05} & 566 & 398\tabularnewline
\hline 
\hline 
\multicolumn{11}{c}{Four-stage Gauss--Legendre scheme with sixth-order finite differences}\tabularnewline
\hline 
\multirow{1}{*}{$\nicefrac{1}{32}$} & 29,791 & \multirow{1}{*}{560,263} & \multirow{1}{*}{5.64} & \multirow{1}{*}{3.56} & \multirow{1}{*}{\textbf{4.77}} & \multirow{1}{*}{\textbf{16}} & \multirow{1}{*}{1.0} & \multirow{1}{*}{\textbf{0.59}} & \multirow{1}{*}{6.57} & \multirow{1}{*}{144}\tabularnewline
\multirow{1}{*}{$\nicefrac{1}{64}$} & 250,047 & \multirow{1}{*}{4,727,079} & \multirow{1}{*}{5.38} & \multirow{1}{*}{35.7} & \multirow{1}{*}{\textbf{57.0}} & \multirow{1}{*}{\textbf{31}} & \multirow{1}{*}{1.0} & \multirow{1}{*}{\textbf{2.77}} & \multirow{1}{*}{95.4} & \multirow{1}{*}{257}\tabularnewline
\multirow{1}{*}{$\nicefrac{1}{128}$} & 2,048,383 & \multirow{1}{*}{38,822,503} & \multirow{1}{*}{5.37} & \multirow{1}{*}{343} & \multirow{1}{*}{\textbf{805}} & \multirow{1}{*}{\textbf{75}} & \multirow{1}{*}{1.0} & \multirow{1}{*}{$-$} & \multirow{1}{*}{$-$} & \multirow{1}{*}{$-$}\tabularnewline
\hline 
\end{tabular}
\end{table}

\section{\label{sec:Conclusion}Conclusion and future work}

In this work, we introduced a software package, \emph{HIFIR}, for
preconditioning GMRES and FGMRES for solving unsymmetric sparse linear
systems. Unlike previous software packages, HIFIR is designed to solve
singular and near-singular (aka ill-conditioned) systems, including
finding least-squares solutions for consistent singular systems, null-space
vectors of singular matrices, and pseudoinverse solutions for inconsistent
systems. This unique feature is backed by a new theory of $\epsilon$-accurate
AGI, and a new algorithm that combines multilevel incomplete LU factorization
with an RRQR on the final Schur complement. Compared to its predecessor
HILUCSI, \emph{HIFIR} also introduces an algorithmic innovation, namely
IBRR, which improves the robustness and significantly reduces the
size of the final Schur complement for some systems. HIFIR was implemented
in C++, with user-friendly high-level interfaces for MATLAB and Python
in \textsf{hifir4m} and \textsf{hifir4py}, respectively. We have released
them as open-source software. We described the software design of
HIFIR in terms of its efficient data structures and its template-based
generic programming interfaces for mixed-precision real and complex
values. We also demonstrated the effectiveness of HIFIR for ill-conditioned
or singular systems arising from several applications, including the
Helmholtz equation, linear elasticity, incompressible Navier--Stokes
equations, and advection-diffusion equation. As presented in this
work, HIFIR was serial. However, it can be used as the computational
kernel in a domain-decomposition preconditioner \cite{smith1996domain}
by factorizing the diagonal blocks within each processor. In addition,
we are presently developing a multi-threaded implementation with the
option of applying the multilevel solver and multilevel matrix-vector
multiplication on GPUs, which we plan to release in the future.
\begin{acks}
Computational results were obtained using the Seawulf cluster at the
Institute for Advanced Computational Science of Stony Brook University,
which was partially funded by the Empire State Development grant NYS
\#28451.
\end{acks}
\bibliographystyle{ACM-Reference-Format}
\bibliography{references}

\appendix

\section{\label{sec:Flexible-array-based}Flexible array-based sparse matrix
data structure}

We describe a flexible, three-tiered, array-based (instead of pointer-based)
sparse matrix data structure to support fan-in updates, deferring,
and rook pivoting, respectively. We focus on the column-oriented version
for $\boldsymbol{L}$; the data structure for $\boldsymbol{U}$ uses
a corresponding row-oriented version.

\subsection{Partially augmented CSC for fan-in updates}

We first briefly describe the baseline data structure as in \cite{li2003crout},
which extends the standard CSC format. Recall that the CSC format
has the following three arrays:
\begin{itemize}
\item \textsf{val}: a floating-point array of size equal to the total number
of nonzeros, with nonzeros stored column by column;
\item \textsf{row\_ind}: an integer array of size equal to that \textsf{of
row\_ind}, storing the row indices in each column;
\item \textsf{col\_start}: an integer array of size $n+1$, storing the
starting index of each column in \textsf{row\_ind}, where $\texttt{col\_start}\left(k\right)$
stores the starting position for the $k$th column in \textsf{row\_ind}.
\end{itemize}
To support fan-in update in Algorithm~\ref{alg:ilu_factor}, we maintain
two additional arrays \textsf{Lstart} and \textsf{Llist} as in \cite{li2003crout},
where \textsf{Lstart} is a size-$n$ integer array storing the first
entry in each column of $\boldsymbol{L}$ whose row index is no smaller
than $k$ at the $k$th step, and \textsf{Llist} is a size-$n$ array-based
linked list storing the index of each column that has a nonzero entry
in $\boldsymbol{\ell}_{k}^{T}$. The combination of \textsf{Lstart}
and \textsf{Llist} allows efficient access of the $k$th row in $\boldsymbol{L}$.
We refer to this baseline data structure as \emph{partially augmenter
CSC} (\emph{PACSC}). The data structure for $\boldsymbol{U}$ uses
the counterpart \emph{PACSR}, which partially augments CSR.

\subsection{Partially augmented CSC with gaps}

\begin{figure}
\includegraphics[width=1\columnwidth]{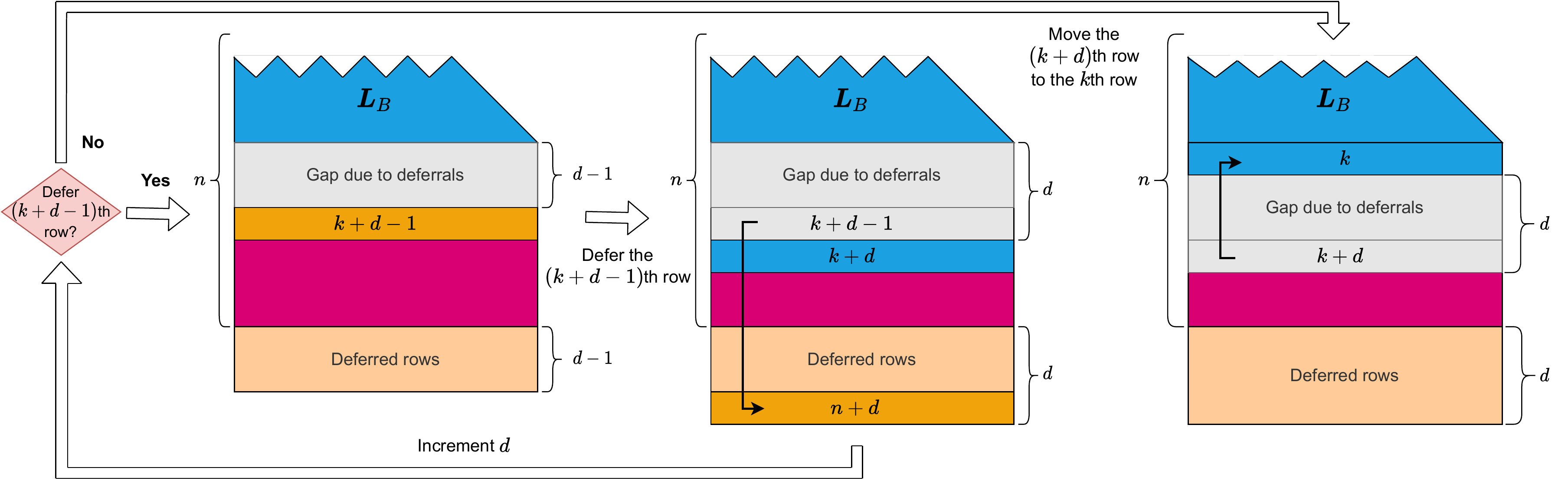}\caption{\label{fig:def-bi-idx}Illustration of PACSC-G (partially augmented
CSC with gap). The left panel shows a gap of $d-1$ before step $k$
of fan-in ILU; the middle panel defers the $k$th row stored in the
$(k+d-1)$th position to the end of $\boldsymbol{L}_{B}$; the right
panel moves the $k$th row from the $(k+d)$th position to the end
of $\boldsymbol{L}_{B}$ for computation.}
\end{figure}

PACSC is memory efficient, but it does not support permutations, such
as static or dynamic deferring. To support static and dynamic deferring,
we extend the data structure to allow a ``gap'' of size $d$ between
$\boldsymbol{L}_{B}$ and $\boldsymbol{L}_{E}$ if there have been
$d$ deferrals. We maintain the gap as follows. At the $k$th step,
suppose there have been $d-1$ deferrals in the preceding steps. Before
performing computation on the $k$th row, we first move $\boldsymbol{\ell}_{k+d-1}^{T}$
to $\boldsymbol{\ell}_{k}^{T}$ in $\boldsymbol{L}$ to eliminate
the gap. If row $k$ needs to be deferred, we directly move the row
from $\boldsymbol{\ell}_{k+d-1}^{T}$ to row $\boldsymbol{\ell}_{n+d}^{T}$
in $\boldsymbol{L}$, which increases the gap to $d$. Figure~\ref{fig:def-bi-idx}
illustrates these operations. To avoid dynamic expansion of the arrays,
we pre-allocate the CSC storage to allow up to $2n$ rows, and also
enlarge \textsf{Llist} and \textsf{Ulist }to size 2$n$. At the end
of \textbf{ilu\_factorize}, we eliminate the gap by moving $\boldsymbol{L}_{n_{1}+d:n+d}$
to $\boldsymbol{L}_{n_{1}:n}$. We refer to the above data structure
as \emph{PACSC-G}; the data structure for $\boldsymbol{U}$ uses the
corresponding \emph{PACSR}-G. It is worth noting that in Algorithm~\ref{alg:ilu_factor},
row $i$ in $\boldsymbol{L}$ for $i>k$ refers to row $i+d$ in PACSC-G,
and similarly, column $j$ in $\boldsymbol{U}$ for $j>k$ refers
to column $j+d$ in PACSR-G. Note that ILUPACK \cite{bollhofer2011ilupack}
also extended CSC to support deferring, but we could not find its
implementation details in its documentation for comparison.

\subsection{Fully augmented CSC for rook pivoting}

To support row interchanges of $\boldsymbol{L}$ in rook pivoting,
we need to access its $i$th row for $k\leq i\leq n_{1}$. The PACSC
(with or without gap) is insufficient for this purpose. To support
rook pivoting, we need to access any row in $\boldsymbol{L}_{k:n,1:k}$
and any column in $\boldsymbol{U}_{1:k,k:n}$. To support this, we
replace the \textsf{Lstart} and \textsf{Llist }arrays in PACSC with
the following additional arrays, analogous to those in the CSR (Compressed
Sparse Row) format:
\begin{itemize}
\item \textsf{val\_pos}: an integer array of size equal to the total number
of nonzeros, storing the indices of the nonzeros in the \textsf{val}
array in the underlying CSC format;
\item \textsf{inv\_val\_pos}: an integer array storing the inverse mapping
of \textsf{val\_pos}.
\item \textsf{col\_ind}: an integer array of size equal to that of \textsf{val\_pos},
storing the column indices within \textsf{val\_pos} in each row;
\item \textsf{row\_start}: an integer array of size $n$, storing the starting
index of each row in \textsf{col\_ind};
\item \textsf{row\_next}: an integer array of the same size as \textsf{col\_ind},
storing the index in \textsf{col\_ind} for the next nonzero in the
same row;
\item \textsf{row\_end}: an integer array of size $n$, storing the last
index of each row in \textsf{col\_ind.}
\end{itemize}
We refer to this data structure as the \emph{fully augmented CSC}
(or \emph{FACSC}). Here, \textsf{val\_pos }and \textsf{inv\_val\_pos}
play the same role as \textsf{val} in CSR; we introduced them to avoid
duplicating the numerical values. \textsf{col\_ind }and \textsf{row\_start}
play the same roles as their respective counterparts in CSR. The arrays
\textsf{row\_next} and \textsf{row\_end }essentially maintain \textsf{Llist}
for all rows in $\boldsymbol{L}_{k:n,1:k}$ instead of for just the
$k$th row. When interchanging rows $i$ and $r$ in $\boldsymbol{L}$
during pivoting, besides updating \textsf{val} and \textsf{row\_ind}
in CSC, we need to swap $\texttt{row\_start}\left(i\right)$ and $\texttt{row\_end}\left(i\right)$
with $\texttt{row\_start}\left(r\right)$ and $\texttt{row\_end}\left(r\right)$,
respectively, while keeping the other arrays intact. The data structure
for $\boldsymbol{U}$ uses \emph{FACSR}, which fully augments CSR.

\section{\label{sec:Time-complexity-rook-pivoting}Time complexity of inverse-based
rook pivoting}

We analyze the time complexity for partial row interchanges in $\boldsymbol{L}$;
the analysis for column interchanges in $\boldsymbol{U}$ is similar.
In Algorithm~\ref{alg:rook_pivot}, there are three key steps: 1)
computing the $k$th column vector $\hat{\boldsymbol{\ell}}$ (line~\ref{line:pivot:ell}),
2) finding a potential pivot $r$ (line~\ref{line:pivot:pvtell}),
and 3) interchanging the $k$th and $r$th rows in $\boldsymbol{L}$
(line~\ref{line:pivot:swapell}). The total number of floating-point
operations in computing $\hat{\boldsymbol{\ell}}$ is bounded by $\mathcal{O}\left(\text{nnz}\left(\hat{\boldsymbol{a}}_{p_{k}}\right)+\underset{i\in\text{nnz}\left(\boldsymbol{u}_{k}\right)}{\sum}\text{nnz}\left(\boldsymbol{L}_{k+1:n,i}\right)\right)$.
Under the assumption that the averaged number of nonzeros per row
and column in the input is a constant (say, bounded by $C$), the
scalability-oriented dropping ensures that this cost is also bounded
by a constant (proportional to $C^{2}$). The inverse-based constraint
in line~\ref{line:pivot:pvtell} also introduces extra cost in estimating
the inverse norm. A brute-force implementation leads to $\mathcal{O}\left(\underset{i\in\text{nnz}\left(\hat{\boldsymbol{\ell}}\right)}{\sum}\text{nnz}\left(\boldsymbol{\ell}_{i}^{T}\right)\right)$
cost in the worst case, which can be improved to $\mathcal{O}\left(\text{nnz}\left(\hat{\boldsymbol{\ell}}\right)\log\left(\text{nnz}\left(\hat{\boldsymbol{\ell}}\right)\right)+\text{nnz}\left(\boldsymbol{\ell}_{r}^{T}\right)\right)$
if we sort $\hat{\boldsymbol{\ell}}$ first and then estimate the
inverse norm only once. Finally, with the FACSC data structure, interchanging
the $k$th row and the $r$th row in $\boldsymbol{L}$ can be done
in $\mathcal{O}(1)$ (or more precisely, $\mathcal{O}(C^{2})$) operations.
Therefore, the time complexity is bounded by constant per row interchange
in rook pivoting. Since we limit the maximum number of row and column
interchanges in rook pivoting by a constant, the total operations
in rook pivoting is no greater than that of the fan-in update in ILU,
and hence rook pivoting does not increase the overall time complexity.
\end{document}